\documentclass{amsart}
\usepackage{amsmath}
\usepackage{amssymb}
\usepackage{amsfonts}

\theoremstyle{plain}

\numberwithin{equation}{section}

\begin{document}
\title[Annuli]{Complex algebraic curves. Annuli}
\author{Maciej Borodzik}
\address{Institute of Mathematics, University of Warsaw, ul. Banacha 2,
02-097 Warsaw, Poland}
\email{mcboro@mimuw.edu.pl}
\author{Henryk \.{Z}o\l \c{a}dek}
\email{zoladek@mimuw.edu.pl}
\date{February 1, 2007}
\thanks{Supported by Polish KBN Grant No 1 P03A 015 29}

\begin{abstract}
We give complete classification of algebraic curves in $\mathbb{C}^{2}$
which are homeomorphic with $\mathbb{C}^{\ast }$ and which satisfy certain
natural condition about codimensions of its singularities. In the proof we
use the method developed in [BZI]. It relies on estimation of certain
invariants of the curve, the so-called numbers of double points hidden at
singularities and at infinity. The sum of these invariants is given by the
Poincar\'{e}--Hopf formula applied to a suitable vector field.
\end{abstract}

\maketitle

\section{The result}

By a plane algebraic annulus we mean a complex reduced algebraic curve $%
\mathcal{C}\subset \mathbb{C}^{2}$ which is a topological embedding of $%
\mathbb{C}^{\ast }=\mathbb{C}\setminus 0.$ Therefore $\mathcal{C}$ has two
places at infinity and can have only cuspidal finite singularities, i.e.
with one local component. Any such curve can be defined by an algebraic
equation $f(x,y)=0,$ but we prefer its parametric definition 
\begin{equation}
x=\varphi (t),\;\;\;y=\psi (t),
\end{equation}%
where $\varphi ,\psi $ are Laurent polynomials and the map 
\begin{equation*}
\xi =(\varphi ,\psi ):\mathbb{C}^{\ast }\rightarrow \mathbb{C}^{2}
\end{equation*}%
defined by (1.1) is one-to-one.

The aim of this work is to classify the algebraic annuli up to \textit{%
equivalence} defined by:

\begin{itemize}
\item polynomial diffeomorphisms of the plane,

\item change of parametrization.
\end{itemize}

Recall that by the Jung--van der Kulk theorem (see [AbMo]) any polynomial
automorphism of $\mathbb{C}^{2}$ (so-called Cremona transformation) is a
composition of a linear map and of \textit{elementary transformations} $%
(x+P(y),y),$ $(x,y+Q(x)).$ The parameter $t$ can be changed to $\lambda t$
or to $\lambda /t.$

In the following theorem we present a list of embedded annuli which satisfy
so-called regularity condition. Roughly speaking, this condition means that
some Puiseux coefficients in Puiseux expansions of local branches of $%
\mathcal{C}$ at the singular points form regular sequences, when treated as
functions on finite dimensional spaces of annuli with fixed asymptotic at
infinity. The regularity condition is defined in the next section and is
studied in our subsequent paper [BZIII].

\bigskip

\textbf{Main Theorem.} \emph{Any algebraic embedding of }$\mathbb{C}^{\ast }$%
\emph{\ into }$\mathbb{C}^{2}$\emph{\ satisfying the regularity condition is
equivalent to one from the below list of pairwise non-equivalent curves (19
series and 4 exceptional cases):}

\emph{(a) }$x=t^{m},$\emph{\ }$y=t^{n}+b_{1}t^{-m}+b_{2}t^{-2m}+\ldots
+b_{k}t^{-km},$\emph{\ }$\gcd (m,n)=1,$\emph{\ }$k=0,1,\ldots ,$\emph{\ }$%
b_{j}\in \mathbb{C},$\emph{\ (}$b_{k}=1$\emph{\ if }$k>0);$

\emph{(b) }$x=t(t-1),$\emph{\ }$y=R_{k,m}(\frac{1}{t}),$\emph{\ }$%
k=1,2,\ldots ,$\emph{\ }$m=0,1,\ldots ,$ $(k,m)\not=(1,0),(2,0),$ $(1,1),$ 
\emph{and } $R_{k,m}$ \emph{are Laurent polynomials defined via }$%
R_{0,m}(u)=(\frac{1}{u}-\frac{1}{2})^{2m+1}$, $%
R_{k+1,m}(u)=[R_{k,m}(u)-R_{k,m}(1)]u^{2}/(u-1);$

\emph{(c) }$x=t^{mn}(t-1),$\emph{\ }$y=S_{k}(\frac{1}{t}),$\emph{\ }$%
k=1,2,\ldots ,$\emph{\ }$n=2,3,\ldots ,$ $mn\geq 2,$ \emph{and }$S_{k}$\emph{%
\ are defined via }$S_{0}(u)=u^{n},$\emph{\ }$%
S_{k+1}(u)=[S_{k}(u)-S_{k}(1)]u^{mn+1}/(u-1);$

\emph{(d) }$x=t^{mn-1}(t-1),$\emph{\ }$y=T_{k}(\frac{1}{t}),$\emph{\ }$%
k=1,2,\ldots ,$\emph{\ }$n=2,3,\ldots ,$ $mn\geq 3,$ \emph{and }$T_{k}$\emph{%
\ are defined via }$T_{0}(u)=u^{n},$\emph{\ }$%
T_{k+1}(u)=[T_{k}(u)-T_{k}(1)]u^{mn}/(u-1);$

\emph{(e) }$x=t^{mn}(t-1),$\emph{\ }$y=U_{k}(\frac{1}{t}),$\emph{\ }$%
k=1,2,\ldots ,$\emph{\ }$n=2,3,\ldots ,$ $mn\geq 2,$ \emph{\ and }$%
U_{0}(u)=u^{-n},$\emph{\ }$U_{k+1}(u)=[U_{k}(u)-U_{k}(1)]u^{mn+1}/(u-1);$

\emph{(f) }$x=t^{mn-1}(t-1),$\emph{\ }$y=V_{k}(\frac{1}{t}),$\emph{\ }$%
k=1,2,\ldots ,\emph{\ }n=2,3,\ldots ,mn\geq 4,$\emph{\ and }$%
V_{0}(u)=u^{-n}, $\emph{\ }$V_{k+1}(u)=[V_{k}(u)-V_{k}(1)]u^{mn}/(u-1);$

\emph{(g) }$x=t^{2}(t-1),$\emph{\ }$y=W_{k}(\frac{1}{t}),$\emph{\ }$%
k=1,2,\ldots ,$\emph{\ and }$W_{1}(u)=3u-u^{2},$\emph{\ }$%
W_{k+1}(u)=[W_{k}(u)-W_{k}(1)]u^{3}/(u-1);$

\emph{(h) }$x=t^{3}(t-1),$\emph{\ }$y=X_{k}(\frac{1}{t}),$\emph{\ }$%
k=1,2,\ldots ,$\emph{\ and }$X_{1}(u)=2u^{2}-u^{3},$\emph{\ }$%
X_{k+1}(u)=[X_{k}(u)-X_{k}(1)]u^{4}/(u-1);$

\emph{(i) }$x=t^{3}(t-1),$\emph{\ }$y=Y_{k}(\frac{1}{t}),$\emph{\ }$%
k=1,2,\ldots ,$\emph{\ and }$Y_{1}(u)=2u^{2}+u^{3},$\emph{\ }$%
Y_{k+1}(u)=[Y_{k}(u)-Y_{k}(1)]u^{4}/(u-1);$

\emph{(j) }$x=Z_{m,n}(t),$\emph{\ }$y=t+\frac{1}{t},\emph{\ }0\leq m\leq n,%
\emph{\ }\left( m,n\right) \neq (0,0),$\emph{\ and the polynomials }$Z_{m,n}$%
\emph{\ are defined by }$Z_{m,n}(t)-Z_{m,n}(\frac{1}{t}%
)=(t-1)^{2m+1}(t+1)^{2n+1}t^{-m-n-1};$

\emph{(k) }$x=(t-1)^{3}t^{-2},$\emph{\ }$y=x^{k}\cdot (t-1)(t-4)t^{-1},$%
\emph{\ }$k=1,2,\ldots ;$

\emph{(l) }$x=(t-1)^{m}t^{-pn},$\emph{\ }$y=(t-1)^{k}t^{-pl},$\emph{\ }$%
ml-nk=1,$\emph{\ }$p=1,2,\ldots ;$

\emph{(m) }$x=(t-1)^{pm}t^{-n},$\emph{\ }$y=(t-1)^{pk}t^{-l},$\emph{\ }$%
ml-nk=1,$\emph{\ }$p=2,3,\ldots ;$

\emph{(n) }$x=(t-1)^{2m}t^{-2n},$\emph{\ }$y=(t-1)^{2k}t^{-2l},$\emph{\ }$%
ml-nk=1;$

\emph{(o) }$x=y^{n}\cdot (t-1)^{2m}(t+1)t^{-m},$\emph{\ }$%
y=(t-1)^{4m}t^{1-2m},$\emph{\ }$m=1,2,\ldots ,$ $n=0,1,\ldots ;$

\emph{(p) }$x=(t-1)^{4}t^{-3},$\emph{\ }$y=x^{k}\cdot (t-1)^{2}(t-3)t^{-2},$%
\emph{\ }$k=0,1,\ldots ;$

\emph{(q) }$x=y^{n}\cdot (t-1)^{2m-1}(t+1)t^{-m},$\emph{\ }$%
y=(t-1)^{4m-2}t^{1-2m},$\emph{\ }$m=2,3,\ldots ,$\emph{\ }$n=0,1,\ldots ;$

\emph{(r) }$x=y^{n}\cdot (t-1)^{3}(t+e^{\pi i/3})t^{-2},$\emph{\ }$%
y=(t-1)^{6}t^{-3},$\emph{\ }$n=0,1,\ldots ;$

\emph{(s)} $x=t^{2n}(t^{2}+\sqrt{2}t+1),$ $y=t^{-2n-4}(t^{2}-\sqrt{2}t+1),$ $%
n=1,2,\ldots ;$

\emph{(t) }$x=(t^{2}+t+\frac{2}{3})t^{4},$\emph{\ }$y=(t^{2}-t+\frac{1}{3}%
)t^{-8};$

\emph{(u) }$x=(t-1)^{2}(t+2)t^{-1},$\emph{\ }$y=(t-1)^{4}(t+\frac{1}{2}%
)t^{-2};$

\emph{(v) }$x=(t-1)^{2}(t+4+2\sqrt{5})t^{-1},$\emph{\ }$y=(t-1)^{4}\left( t+%
\frac{1}{4}\left( 11+5\sqrt{5}\right) \right) t^{-2};$

\emph{(w) }$x=(t-1)^{2}(t+2)t^{-1},$\emph{\ }$y=(t-1)^{2}(t+\frac{1}{2}%
)t^{-2}.$

\bigskip

\textbf{Commentary}. Here we present singularities of the curves listed in
Main Theorem. We expose the essential terms in expansions of these curves at
the singular points in the affine parts of the curves as well as at the
infinity.

(a) The curve is smooth, i.e. in the affine part. As $t\rightarrow \infty $
we have $x\sim t^{m}$ and $y\sim x^{n/m}+$(integer powers of $x).$ As $%
t\rightarrow 0$ we have $x\sim t^{m}$ and $y\sim x^{n/m}+$(integer powers of 
$x).$ If $n>m$ the curve has only one point at infinity (in $\mathbb{CP}^{2}$%
), otherwise there are two such points. See also Lemma 3.2.

(b) The curve has one singular point of the type \textbf{A}$_{2m},$ i.e. $%
Y^{2}=X^{2m+1},$ at $t=\frac{1}{2}.$ As $t\rightarrow \infty $ we have $%
x\sim t^{2}$ and $y\sim c_{1}x^{-1}+c_{2}x^{m-k+1/2}$ (here $m-k+1/2$ can be 
$>-1).$ As $t\rightarrow 0$ we have $x\sim t$ and $y\sim t^{-k}.$ See also
Lemma 3.10.

(c) The curve is smooth. As $t\rightarrow \infty $ we have $x\sim t^{mn+1}$
and $y\sim c_{1}x^{-1}+c_{2}x^{-k-n/(mn+1)}.$ As $t\rightarrow 0$ we have $%
x\sim t^{mn}$ and $y\sim x^{-k-1/m}(1+cx^{1/mn}).$ See also Lemma 3.11.

(d) The curve is smooth. As $t\rightarrow \infty $ we have $x\sim t^{mn}$
and $y\sim c_{1}x^{-1}+c_{2}x^{-k-1/m}(1+c_{3}x^{-1/mn}).$ As $t\rightarrow
0 $ we have $x\sim t^{mn-1}$ and $y\sim x^{-k-n/(mn-1)}.$ See also Lemma
3.11.

(e) The curve is smooth. As $t\rightarrow \infty $ we have $x\sim t^{mn+1}$
and $y\sim c_{1}x^{-1}+c_{2}x^{-k+n/(mn+1)}.$ As $t\rightarrow 0$ we have $%
x\sim t^{mn}$ and $y\sim c_{1}x^{-1}+c_{2}x^{-k+1/m}(1+c_{3}x^{1/mn}).$ See
also Lemma 3.14.

(f) The curve is smooth. As $t\rightarrow \infty $ we have $x\sim t^{mn}$
and $y\sim c_{1}x^{-1}+c_{2}x^{-k+1/m}(1+c_{3}x^{-1/mn}).$ As $t\rightarrow
0 $ we have $x\sim t^{mn-1}$ and $y\sim c_{1}x^{-1}+c_{2}x^{-k+n/(mn-1)}.$
See also Lemma 3.14.

(g) There is the cusp singularity \textbf{A}$_{2}$ at $t=2/3.$ As $%
t\rightarrow \infty $ we have $x\sim t^{3}$ and $y\sim cx^{-1}+x^{-k-1/3}.$
As $t\rightarrow 0$ we have $x\sim t^{2}$ and $y\sim
c_{1}x^{-1}+c_{2}x^{-k}(1+c_{3}x^{1/2}).$ See also Lemma 3.15.

(h) There is the cusp singularity \textbf{A}$_{2}$ at $t=3/4$. As $%
t\rightarrow \infty $ we have $x\sim t^{4}$ and $y\sim
c_{1}x^{-1}+c_{2}x^{-k+1/2}(1+c_{3}x^{-1/4}).$ As $t\rightarrow 0$ we have $%
x\sim t^{3}$ and $y\sim c_{1}x^{-1}+c_{2}x^{-k+1/3}).$ See also Lemma 3.15.

(i) The curve is smooth. As $t\rightarrow \infty $ we have $x\sim t^{4}$ and 
$y\sim c_{1}x^{-1}+c_{2}x^{-k+1/2}(1+c_{3}x^{-3/4}),$ i.e. we have a
degeneration. As $t\rightarrow 0$ we have $x\sim t^{3}$ and $y\sim
c_{1}x^{-1}+c_{2}x^{-k+1/3}).$ See also Lemma 3.15.

(j) As $t\rightarrow \infty $ we have $y\sim t$ and $x\sim y^{m+n+1}$ and as 
$t\rightarrow 0$ we have $y\sim t^{-1}$ and $x\sim y^{-1}$ (smoothness). It
has two singular points: at $t=1$ and at $t=-1$ of the type \textbf{A}$_{2m}$
(i.e. $(t^{2},t^{2m+1}))$ and of the type \textbf{A}$_{2n}$ respectively.
See also Lemma 4.18.

(k) The curve has two singular points: at $t=1$ with $x\sim (t-1)^{3},$ $%
y\sim x^{k+1/3}$ and the cusp \textbf{A}$_{2}$ at $t=-2.$ As $t\rightarrow
\infty $ we have $x\sim t$ and $y\sim x^{k+1}.$ As $t\rightarrow 0$ we have $%
x\sim t^{-2}$ and $y\sim x^{-k-1/2}.$ See also Lemma 5.14.

(l) The curve has singular point at $t=1:$ $x\sim (t-1)^{m}$ and $y\sim
x^{k/m}$. As $t\rightarrow \infty $ we have $x\sim t^{m-pn}$ and $y\sim
t^{k-pl}.$ As $t\rightarrow 0$ we have $x\sim t^{-pn}$ and $y\sim
x^{l/n}(1+cx^{-1/pn}).$ See also Lemma 5.15.

(m) The curve has singular point at $t=1:$ $x\sim (t-1)^{pm}$ and $y\sim
x^{k/m}(1+cx^{1/pm})$. As $t\rightarrow \infty $ we have $x\sim t^{pm-n}$
and $y\sim t^{pk-l}.$ As $t\rightarrow 0$ we have $x\sim t^{-n}$ and $y\sim
t^{-l}.$ See also Lemma 5.15.

(n) The curve has singular point at $t=1:$ $x\sim (t-1)^{2m}$ and $y\sim
x^{k/m}(1+cx^{1/2m})$. As $t\rightarrow \infty $ we have $x\sim t^{2(m-n)}$
and $y\sim x^{(k-l)/(m-n)}(1+cx^{-1/2(m-n)}).$ As $t\rightarrow 0$ we have $%
x\sim t^{-2n}$ and $y\sim x^{l/n}(1+cx^{-1/2n}).$ See also Lemma 5.15.

(o) The curve has singular point at $t=1:$ $y\sim (t-1)^{4m}$ and $x\sim
y^{n+1/2}(1+c_{1}y^{1/2m}+c_{2}y^{3/4m})$. As $t\rightarrow \infty $ we have 
$y\sim t^{2m+1}$ and $x\sim y^{n+(m+1)/(2m+1)}.$ As $t\rightarrow 0$ we have 
$y\sim t^{1-2m}$ and $x\sim y^{n+m/(2m-1)}.$ See also Lemma 5.20.

(p) The curve has two singular points: at $t=1$ with $x\sim (t-1)^{4}$ and $%
y\sim x^{k+1/2}(1+cx^{1/4})$ and the cusp \textbf{A}$_{2}$ at $t=-3$. As $%
t\rightarrow \infty $ we have $x\sim t$ and $y\sim x^{k+1}$ and $x\sim
t^{-3} $, $y\sim x^{k+2/3}$ as $t\rightarrow 0.$ See also Lemma 5.14.

(q) The curve has singular point at $t=1$ with $y\sim (t-1)^{2(2m-1)}$ and $%
x\sim y^{n+1/2}(1+cy^{2/2(2m-1)})$. As $t\rightarrow \infty $ we have $y\sim
t^{2m-1}$ and $x\sim y^{n+m/(2m-1)}.$ As $t\rightarrow 0$ we have $y\sim
t^{1-2m}$ and $x\sim -y^{n+m/(2m-1)}.$ Thus the two local branches at
infinity have the same order of the asymptotic but differ in the leading
coefficient. See also Lemma 5.37.

(r) The curve has singular point at $t=1$ with $y\sim (t-1)^{6}$ and $x\sim
y^{n+1/2}(1+cy^{1/6)})$. As $t\rightarrow \infty $ we have $y\sim t^{3}$ and 
$x\sim y^{n+2/3}.$ As $t\rightarrow 0$ we have $y\sim t^{-3}$ and $x\sim
(-e^{i\pi /3})y^{n+2/3}.$ Thus the two local branches at infinity have the
same order of the asymptotic and the leading coefficients, denoted $A$ and $%
B $ respectively, satisfy $A^{3}=B^{3}$. It is the whole degeneration at
infinity. See also Lemma 5.37.

(s) The curve is smooth. As $t\rightarrow \infty $ we have $x\sim t^{2(n+1)}$
and $y\sim c_{1}x^{-1}+c_{2}x^{-1-2/(n+1)}(1+c_{3}x^{-1/2(n+1)})$. As $%
t\rightarrow 0$ we have $x\sim t^{2n}$ and $y\sim x^{-(n+2)/n}(1+cx^{1/2n}).$
See also Lemma 3.20.

(t) The curve is smooth. As $t\rightarrow \infty $ we have $x\sim t^{6}$ and 
$y\sim x^{-1}+c_{1}x^{-3/2}(1+c_{2}x^{-1/6})$. As $t\rightarrow 0$ we have $%
x\sim t^{4}$ and $y\sim c_{1}x^{-2}+c_{2}x^{-3/2}+c_{3}x^{-5/4}.$ See also
Lemma 3.27.

(u) The curve has singularity of the type \textbf{A}$_{8}$ at $t=1$. As $%
t\rightarrow \infty $ we have $x\sim t^{2}$ and $y\sim x^{3/2}$. As $%
t\rightarrow 0$ we have $x\sim t^{-1}$ and $y\sim x^{2}$ (smoothness). See
also Lemma 5.8.

(v) The curve has two singularities of the type \textbf{A}$_{4}:$ at $t=1$
and at $t=\frac{2}{3}(\sqrt{5}-2)$. As $t\rightarrow \infty $ we have $x\sim
t^{2}$ and $y\sim x^{3/2}$. As $t\rightarrow 0$ we have $x\sim t^{-1}$ and $%
y\sim x^{2}$ (smoothness). See also Lemma 5.8.

(w) The curve has three cusps \textbf{A}$_{2}$ and is smooth at $t=\infty $ (%
$y\sim t,$ $x\sim y^{2})$ and at $t=0$ ($x\sim t^{-1},$ $y\sim x^{2}).$ See
also Lemma 6.1.

$\smallskip $

It follows that exactly in the cases (a), (c), (d), (e), (f), (i), (s) and
(t) the $\mathbb{C}^{\ast }$--embedding is smooth.

We have checked that any curve from the above list can be reduced to a
straight line by means of a birational change of $\mathbb{CP}^{2}.$ (The
same holds for affine rational curves with one self-intersection, which were
classified in [BZI].) We do not present these changes; the reader can easily
do it case by case. This confirms the conjecture that any rational curve in $%
\mathbb{CP}^{2}$ can be straightened via a birational automorphism (see
[FlZa]).

In contrast to the list given in [BZI] the classification from Main Theorem
contains moduli. These moduli $b_{2},\ldots ,b_{l}$ appear only in the case
(a). Our method does not explanation this phenomenon in a satisfactory way.

\bigskip

Recall that Main Theorem does not yet solve the problem of classification of
(topological) embeddings of $\mathbb{C}^{\ast }$ into $\mathbb{C}^{2}.$ It
assumes some bound on codimensions of singularities of (topological)
immersions of $\mathbb{C}^{\ast }$ into $\mathbb{C}^{2}$ which are stated in
Conjecture 2.40 below. In this sense Main Theorem is an analogue of the main
result of our previous paper [BZI], where a classification of
(topologically) immersed lines $\mathbb{C}$ into $\mathbb{C}^{2}$ with one
self-intersection point is given (21 cases with 16 series and 5 exceptional
cases) under an analogous assumption about codimensions. In our forthcoming
paper [BZIII] we prove some results about the codimensions (see also Remark
2.43 below). The bounds obtained are not optimal, the discrepancy between
this bound and the dimension of (topologically) immersed annuli is $\geq 4.$
In principle it is possible to complete the proof of the classification from
Main Theorem, for this one has to analyze a lot more cases that in the below
proof. We have done this analysis (without publication) for annuli of the
type $\binom{-}{-}$ (see (2.35)) and no new cases have appeared.

We recall that any simply connected curve either is rectifiable $x=t,$ $y=0$
(the Abhyankar--Moh--Suzuki theorem [AbMo], [Su1]) or is equivalent to the
quasi-homogeneous curve $x=t^{k},$ $y=t^{l}$ (the Zaidenberg--Lin theorem
[ZaLi]).

There are not many results about curves homeomorphic to an annulus. W.
Neumann [Ne] proved that if such curve $f=0$ is smooth and typical in the
family $f=\lambda ,$ then it is equivalent to the case (a) of Main Theorem.
L. Rudolph in [Ru] gives the example $x=t^{2}+2t^{-2},$ $y=2t+t^{-2}$ of a
projective rational cuspidal curve with three cusps; it is the case (w) of
Main Theorem.

S. Kaliman [Ka] classified all smooth embeddings of $\mathbb{C}^{\ast }$
into $\mathbb{C}$ such that the corresponding polynomial $F$ has rational
level curves. These are $F(x,y)=$\linebreak $\left[ \chi ^{mn+1}-(\chi
^{n}+x)^{m}\right] /x^{m}=0$ and $\left[ \chi ^{mn-1}-(\chi ^{n}+x)m\right]
/x^{m}=0,$ where $\chi =x^{m}y+a_{m-1}x^{m-1}+\ldots +a_{1}x+1$ are such
that the above functions are polynomials. Moreover, $m\geq 2,$ $n\geq 1$ and 
$(m,n)\neq (2,1)$ in the case of second curve. Later we shall see that these
curves correspond to $x=t^{mn}(t-1),$ $y=U_{m}(\frac{1}{t})$ (from the
series (e)) and $x=t^{mn-1}(t-1),$ $y=V_{m}(\frac{1}{t})$ (from the series
(f)) respectively.

The series (s) was firstly found by P. Cassou-Nogu\`{e}s (we owe this
information to M. Koras). Also M. Koras and P. Russell proved (but have not
published yet) that any smooth annulus can be reduced to one of the curves
found by M. Zaidenberg and V. Lin. The problem of classification of annuli
is raised also in the work [NeNo] of W. Neumann and P. Norbury.

\medskip

Among other methods in the study of affine algebraic curves it is worth to
mention the knot invariants (so-called splice diagrams introduced in [EiNe])
used by Neumann and Rudolph (see [NeRu]). The splice diagrams are related
with the dual graphs of the resolution of singularities and of indeterminacy
(of the polynomials defining the curves) at infinity (see [ABCN).

There are some works devoted to study projective curves, see [FlZa], [MaSa],
[Su2], [Or], [OZ1], [OZ2], [Yo], [ZaOr] for example. We do not consider
projective curves (only affine ones), because our method ceases to be
effective in the projective case.\medskip

The method used in this paper was developed in [BZI]. It relies on estimates
of numerical invariants of local singularities, like the Milnor number (or,
better, the number of double points hidden at a singularity), in terms of
suitably defined codimensions of the singularities. The sum of the numbers
of hidden double points is calculated by means of the Poincar\'{e}--Hopf
formula and the sum of codimensions is estimated by the dimension of the
space of parametric rational curves with fixed asymptotic behaviour at
infinity. Such estimates allow to reduce the set of curves, which are
candidates for $\mathbb{C}^{\ast }$-embeddings. There remain several classes
of curves which are studied separately. The details of the method are given
in the next section.

The very proof of Main Theorem is given in Sections 3, 4, 5 and 6, each
devoted to one type of curves.

\bigskip

\section{Estimates for annuli}

\textbf{2.I. The Poincar\'{e}--Hopf formula.} Let $\mathcal{C}=\left\{
f=0\right\} $ be a reduced curve in $\mathbb{C}^{2}.$ The Hamiltonian vector
field $X_{f}=f_{y}\frac{\partial }{\partial x}-f_{x}\frac{\partial }{%
\partial y}$ is tangent to $\mathcal{C}$.

Suppose that $z$ is a singular point of $\mathcal{C}$. Consider the local
normalization $N_{z}:\widetilde{A}\rightarrow (\mathcal{C},z),$ where $%
\widetilde{A}$ is a disjoint union of discs $\widetilde{A}_{j},$ $j=1,\ldots
,k,$ $\widetilde{A}_{j}\simeq \left\{ |z|<1\right\} ,$ such that $%
A_{j}=N_{z}(\widetilde{A}_{j})$ are local irreducible components of $(%
\mathcal{C},z).$ The pull-back $\widetilde{X}=N_{z}^{\ast
}X_{f}=(N_{z})_{\ast }^{-1}X_{f}\circ N_{z}$ of the Hamiltonian vector field
is a vector field on the smooth manifold with isolated singular points $%
p_{j}=N_{z}^{-1}(z)\cap \widetilde{A}_{j},$ $j=1,\ldots ,k.$ Therefore one
can define the indices $i_{p_{j}}\widetilde{X}.$

\bigskip

2.1. \textbf{Definition.} We call the quantity 
\begin{equation*}
\delta _{z}=\frac{1}{2}\sum_{j}i_{p_{j}}\widetilde{X}
\end{equation*}%
the number of double points of $\mathcal{C}$ hidden at $z.$

\bigskip

It is known that (see [BZI]) 
\begin{equation}
2\delta _{z}=\sum_{j}\mu _{z}(A_{j})+2\sum_{i<j}(A_{i}\cdot A_{j})_{z} 
\tag{2.2}
\end{equation}%
and 
\begin{equation}
2\delta _{z}=\mu _{z}(\mathcal{C})+k-1.  \tag{2.3}
\end{equation}%
Here $\mu _{z}(\cdot )$ and $(A_{i}\cdot A_{j})_{z}$ denote the Milnor
number and the intersection index respectively. Therefore $\delta _{z}$
coincides with the standard definition of the number of double points (see
[Mil]). It is the number of double points of a generic perturbation of the
normalization map $N.$ \medskip

Consider now the normalization $N:\widetilde{\mathcal{C}}\rightarrow 
\overline{\mathcal{C}}$ of the closure $\overline{\mathcal{C}}\subset 
\mathbb{CP}^{2}$ of $\mathcal{C}$. The vector field $N^{\ast }X_{f}$ is not
regular, it has poles. Therefore we choose 
\begin{equation*}
\widetilde{X}=h\cdot N^{\ast }X_{f},
\end{equation*}%
where $h:\mathcal{C}\rightarrow \mathbb{R}_{+}$ is a smooth function tending
to zero sufficiently fast near the preimages of the points of $\mathcal{C}$
at infinity. The indices of $\widetilde{X}$ at the preimages of the points
at infinity are well defined.

The \textit{Poincar\'{e}--Hopf formula} states that 
\begin{equation*}
\sum_{t\text{ singular}}i_{t}\widetilde{X}=\chi (\widetilde{\mathcal{C}}), 
\notag
\end{equation*}%
where $\chi (\widetilde{\mathcal{C}})$ denotes the Euler--Poincar\'{e}
characteristic; (it is sometimes called the intristic Euler--Poincar\'{e}
characteristic of $\overline{\mathcal{C}}$).

We are interested in the case when $\widetilde{\mathcal{C}}=\mathbb{CP}^{1}$
and $N^{-1}(\mathcal{C})=\mathbb{CP}^{1}\setminus \left\{
(0:1),(1:0)\right\} =\mathbb{C}^{\ast }.$ The normalization map $N|_{\mathbb{%
C}^{\ast }}$ coincides with the parametrization $t\rightarrow \xi (t).$
\medskip

The number of double points hidden at a cuspidal singularity is expressed
via its Puiseux expansion. Assuming that the curve is locally given by 
\begin{equation}
x=\tau ^{n},\;\;\;y=C_{1}\tau +C_{2}\tau ^{2}+\ldots  \tag{2.4}
\end{equation}
(it is the so-called \textit{standard Puiseux expansion}) we define the 
\textit{topologically arranged Puiseux series}%
\begin{equation}
\begin{array}{c}
y=x^{m_{0}}(D_{0}+\ldots )+x^{m_{1}/n_{1}}(D_{1}+\ldots )+\ldots
+x^{m_{l}/n_{1}\ldots n_{l}}(D_{l}+\ldots ) \\ 
=\tau ^{v_{0}}(D_{0}+\ldots )+\tau ^{v_{1}}(D_{1}+\ldots )+\ldots +\tau
^{v_{l}}(D_{l}+\ldots ).%
\end{array}
\tag{2.5}
\end{equation}%
Here $m_{0}$ is an integer, the \textit{characteristic pairs} $(m_{j},n_{j})$
satisfy $n_{j}>1,$ $\gcd (m_{j},n_{j})=1,$ $n=n_{1}\ldots n_{l}$, $%
v_{0}<v_{1}<\ldots <v_{l}$, the \textit{essential Puiseux coefficients} $%
D_{j}\neq 0$ and the dots in the $j$-th summand mean terms with $%
x^{k/n_{1}\ldots n_{j}}.$ The first summand may be absent.

\bigskip

2.6. \textbf{Proposition} ([Mil]). \emph{We have} 
\begin{eqnarray*}
\mu _{0} &=&2\delta _{0}=\sum_{j=1}^{l}(v_{j}-1)(n_{j}-1)n_{j+1}\ldots n_{l}
\\
&=&\sum (m_{j}n_{j+1}\ldots n_{l}-1)(n_{j}-1)n_{j+1}\ldots n_{l}.
\end{eqnarray*}%
\medskip

The annulus $\mathcal{C}$ has two places at infinity, one corresponding to $%
t=0$ and one corresponding to $t=\infty .$

Let $\mathcal{C}_{\infty }$ denote the branch corresponding to $t\rightarrow
\infty $ with local variable $\tau =t^{-1}:$%
\begin{equation}
\mathcal{C}_{\infty }:x=\tau ^{-p}+\ldots ,\;\;y=\tau ^{-q}+\ldots  \tag{2.7}
\end{equation}%
with the topologically arranged Puiseux expansion 
\begin{equation*}
\mathcal{C}_{\infty }:y=x^{q_{1}/p_{1}}(E_{1}+\ldots
)+x^{q_{2}/p_{1}p_{2}}(E_{2}+\ldots )+\ldots +x^{q_{l\infty }/p_{1}\ldots
p_{l_{\infty }}}(E_{l_{\infty }}+\ldots ),
\end{equation*}%
where $\gcd (q_{j},p_{j})=1$ for the corresponding characteristic pairs.

Let $\mathcal{C}_{0}$ be the second branch:%
\begin{equation}
\mathcal{C}_{0}:x=t^{-r}+\ldots ,\;\;y=t^{-s}+\ldots ,\;\;t\rightarrow 0, 
\tag{2.8}
\end{equation}%
with the Puiseux expansion%
\begin{equation*}
\mathcal{C}_{0}:y=x^{s_{1}/r_{1}}(F_{1}+\ldots
)+x^{s_{2}/r_{1}r_{2}}(F_{2}+\ldots )+\ldots +x^{s_{l_{0}}/r_{1}\ldots
r_{l_{0}}}(F_{l_{0}}+\ldots ),
\end{equation*}%
where $\gcd (s_{j},r_{j})=1.$ The following result is proved in the same way
as Theorem 2.7 in [BZI].

\bigskip

2.9. \textbf{Proposition}. \emph{If }$ps-rq\neq 0,$\emph{\ then }%
\begin{equation*}
i_{\infty }\widetilde{X}=\left\{ 2-\sum_{j=1}^{l_{\infty
}}(q_{j}p_{j+1}\ldots p_{l_{\infty }}-1)(p_{j}-1)p_{j+1}\ldots p_{l_{\infty
}}\right\} -\max (ps,rq)
\end{equation*}%
\emph{and} 
\begin{equation*}
i_{0}\widetilde{X}=\left\{ 2-\sum_{j=1}^{r_{0}}(s_{j}r_{j+1}\ldots
r_{l_{0}}-1)(r_{j}-1)r_{j+1}\ldots r_{l_{0}}\right\} -\max (ps,rq).
\end{equation*}%
\medskip

Denote%
\begin{equation}
p^{\prime }=\gcd (p,q),\;\;\;\;r^{\prime }=\gcd (r,s),  \tag{2.10}
\end{equation}%
in (2.7) and (2.8) they are equal $p_{2}\ldots p_{l_{\infty }}$ and $%
r_{2}\ldots r_{l_{0}}$ respectively. Assuming that the curve is typical,
i.e. that the only singularities are nodal (double) points and that there
are two characteristic pairs at $t=\infty $ ($(q_{1},p_{1})$ and $%
(q-1,p^{\prime }))$ and at $t=0$ ($(s_{1},r_{1})$ and $(s-1,r^{\prime })),$
we get 
\begin{eqnarray*}
i_{0}\widetilde{X}+i_{\infty }\widetilde{X} &=&\left\{ 2-(q-1)(p-p^{\prime
})-(q-2)(p^{\prime }-1)\right\} \\
&&+\left\{ 2-(s-1)((r-r^{\prime })-(s-2)(r^{\prime }-1)\right\} -2\max
(ps,rq) \\
&=&\left\{ 2-\left[ (q-1)(p-1)-(p^{\prime }-1)\right] \right\} +\left\{ 2-%
\left[ (s-1)(r-1)-(r^{\prime }-1)\right] \right\} \\
&&-2\max (ps,rq) \\
&=&2-\left\{ (p+r-1)(q+s-1)-(p^{\prime }+r^{\prime }-1)+\left\vert
ps-rq\right\vert \right\} .
\end{eqnarray*}

It is natural to introduce the \textit{maximal number of double points} $%
\delta _{\max }$ by 
\begin{equation}
2\delta _{\max }=(p+r-1)(q+s-1)-(p^{\prime }+r^{\prime }-1)+\left\vert
ps-rq\right\vert  \tag{2.11}
\end{equation}%
which is the number of finite double points in the typical case (this notion
is valid also when $ps=rq).$ Define also the \textit{number of double points
hidden at} $t=0$, i.e. $\delta _{0},$ by 
\begin{equation}
2\delta _{0}=(2-i_{0}\widetilde{X})-2\delta _{0,\max },\;\;\;2\delta
_{0,\max }=(r-1)(s-1)-(r^{\prime }-1)+\max (ps,rq),  \tag{2.12}
\end{equation}%
the \textit{number of double points hidden at} $t=\infty ,$ $\delta _{\infty
}$ by%
\begin{equation}
2\delta _{\infty }=(2-i_{\infty }\widetilde{X})-2\delta _{\infty ,\max
},\;\;\;2\delta _{\infty ,\max }=(p-1)(q-1)-(p^{\prime }-1)+\max (ps,rq), 
\tag{2.13}
\end{equation}%
and the \textit{number of double points hidden at infinity} 
\begin{equation*}
\delta _{\inf }=\delta _{0}+\delta _{\infty }\;\;\text{when}\;\;ps\neq rq.
\end{equation*}%
The numbers $\delta _{0},\delta _{\infty },$ and $\delta _{\inf }$ control
the degenerations of a given curve at infinity.

\medskip

2.14. \textbf{Remark}. The number $\delta _{\inf }$ should be not confused
with the number of double points hidden at the singular points at infinity
in $\mathbb{CP}^{2}$ of the projective closure $\overline{\mathcal{C}}$ of
the curve $\mathcal{C}$.

\bigskip

2.15. \textbf{Proposition.} \emph{We have }%
\begin{equation}
2\delta _{\inf }+\sum_{P_{j}}2\delta _{p_{j}}=2\delta _{\max },  \tag{2.16}
\end{equation}%
\emph{where the sum runs over finite singular points }$P_{j}=\xi (t_{j})$%
\emph{\ of the curve }$\mathcal{C}=\xi (\mathbb{C}^{\ast }).$\emph{\ It
implies that, for an embedding }$\xi $\emph{\ (with fixed asymptotic as }$%
t\rightarrow 0$\emph{\ and }$t\rightarrow \infty ),$\emph{\ the double
points (for a generic immersion }$\mathbb{C}^{\ast }\rightarrow \mathbb{C}%
^{2})$\emph{\ hide at infinity and at the finite cuspidal singularities.}

\emph{The identity (2.16) holds true also in the case }$ps=rq,$\emph{\ but
with another interpretation of }$\delta _{\inf }$\emph{\ given below.}

\bigskip

In the case $ps-rq=0$ some terms of the Puiseux expansion for the branches $%
\mathcal{C}_{0,\infty }$ may coincide. We have 
\begin{equation}
\begin{array}{c}
\mathcal{C}_{\infty }:x=t^{v\tilde{p}}+\ldots ,\;\;y=G_{1}x^{w/v}+\ldots
+G_{u}x^{(w-u+1)/v}+Ex^{u_{\infty }/v\tilde{p}_{1}}+\ldots \\ 
\mathcal{C}_{0}:x=t^{-v\tilde{r}}+\ldots ,\;\;y=G_{1}x^{w/v}+\ldots
+G_{u}x^{(w-u+1)/v}+Fx^{u_{0}/v\tilde{r}_{1}}+\ldots%
\end{array}
\tag{2.17}
\end{equation}%
Here $u$ terms of the two Puiseux expansions coincide and we assume that it
is maximal such sequence (when taken into account different choices of the
roots $x^{j/v}).$ The terms $Ex^{u_{\infty }/v\tilde{p}_{1}},$ $\gcd
(u_{\infty },\tilde{p}_{1})=1,$ and $Fx^{u_{0}/v\tilde{r}_{1}},$ $\gcd
(u_{0},\tilde{r}_{1})=1,$ are different: either $u_{\infty }/v\tilde{p}%
_{1}\neq u_{0}/v\tilde{r}_{1}$ or $u_{\infty }/v\tilde{p}_{1}=u_{0}/v\tilde{r%
}_{1}$ but $E^{v\tilde{p}_{1}}\neq F^{v\tilde{p}_{1}}.$ Moreover, $v$ is
maximal possible (so that $G_{u}$ might possibly be zero).

Let us arrange topologically the coinciding terms%
\begin{equation*}
\left( \widetilde{G}_{1}x^{w_{1}/v_{1}}+\ldots \right) +\ldots +\left( 
\widetilde{G}_{r}x^{w_{l}/v_{1}\ldots v_{l}}+\ldots \right) ,\;\;\gcd
(w_{j},v_{j})=1.
\end{equation*}

The next result is analogous to Proposition 2.13 in [BZI].

\bigskip

2.18. \textbf{Proposition}. \emph{We have} 
\begin{eqnarray*}
2-i_{0}\widetilde{X}-i_{\infty }\widetilde{X} &=&\sum (q_{j}p_{j+1}\ldots
p_{l_{\infty }}-1)(p_{j}-1)p_{j+1}\ldots p_{l_{\infty }} \\
&&+\sum (s_{j}r_{j+1}\ldots r_{l_{0}}-1)(r_{j}-1)r_{j+1}\ldots r_{l_{0}} \\
&&+2pq\left( \sum w_{j}(v_{j}-1)\left( v_{j+1}\ldots v_{l}\right) ^{2}+\max
\left\{ \frac{u_{\infty }}{\tilde{p}_{1}},\frac{u_{0}}{\tilde{r}_{1}}%
\right\} \right) .
\end{eqnarray*}%
\emph{When we define the number of double points hidden at infinity by}%
\begin{equation}
2\delta _{\inf }=2\delta _{\max }-(2-i_{\infty }\widetilde{X}-i_{0}%
\widetilde{X}),  \tag{2.19}
\end{equation}%
\emph{then the identity (2.16) holds true}.

\bigskip

\textbf{2.II. Bounds for the numbers of hidden double points.} The success
of the paper [BZI] relied upon using very effective estimates for the Milnor
numbers of singularities and for the number of double points hidden at
infinity. Following [BZI] for local cuspidal singularities of the form \emph{%
\ }$x=\tau ^{n},$\emph{\ }$y=C_{1}\tau +C_{2}\tau ^{2}+\ldots ,$ i.e. with
fixed $n,$ we define the \textit{codimension} $\nu $ of the stratum $\mu =$%
const as the number of equations $C_{i}=0$ (vanishing \textit{essential
Puiseux quantities}) appearing in definition of the equisingularity stratum.

\bigskip

2.20. \textbf{Proposition}. ([BZI]) \emph{The Milnor number of such
singularity satisfies}%
\begin{equation}
\mu \leq n\nu .  \tag{2.21}
\end{equation}

\emph{Moreover, when we restrict the class of curves to }%
\begin{equation}
x=\tau ^{n},\;\;y=\tau ^{m}(1+C_{1}\tau +\ldots ),  \tag{2.22}
\end{equation}%
\emph{then }%
\begin{equation}
\mu \leq \mu _{\min }+n^{\prime }\nu ^{\prime },  \tag{2.23}
\end{equation}%
\emph{where }$n^{\prime }=\gcd (m,n),$\emph{\ }$\nu ^{\prime }$\emph{\ is
the codimension of stratum }$\mu =$const\emph{\ and the minimal Milnor
number equals}%
\begin{equation}
\mu _{\min }=(m-1)(n-1)-(n^{\prime }-1).  \tag{2.24}
\end{equation}%
\medskip \bigskip

We complete Proposition 2.20 with presentation of some situations when the
bounds (2.21) and (2.23) become equalities (without straightforward proofs).

\bigskip

2.25. \textbf{Lemma}.\emph{The equality }$\mu =n\nu $\emph{\ holds only in
two cases:}

\emph{(i) when there is only one characteristic pair }$(m,n)$\emph{\ with }$%
m=1$ $\pmod n$\emph{\ (it is always so when }$n=2);$

\emph{(ii) when there are two characteristic pairs }$(m_{1},n_{1}),$\emph{\ }%
$m_{1}=1$ $\pmod{n_{1}}$\emph{\ and }$(m_{1}n^{\prime }+1,n').$

\bigskip

2.26. \textbf{Lemma}. \emph{For the curve (2.22) necessary conditions for
the equality }$\mu =\mu _{\min }+n^{\prime }\nu ^{\prime }$\emph{\ are
following:}

\emph{(i) if there are two characteristic pairs }$(m_{1},n_{1})$\emph{\ and }%
$(m^{\prime },n^{\prime }),$\emph{\ then }$m^{\prime }=1$ $\pmod {n'};$

\emph{(ii) if }$\nu ^{\prime }=1,$\emph{\ then }$n^{\prime }$\emph{\ is even
and }$C_{1}=0;$

\emph{(iii) if }$\nu ^{\prime }=2,$\emph{\ then either }$n^{\prime }=2$\emph{%
\ and }$C_{1}=C_{3}=0,$\emph{\ or }$n^{\prime }=0$ $\pmod 3$\emph{\
and }$C_{1}=C_{2}=0.$\newline
\emph{There are also natural inequalities for suitable essential Puiseux
coefficients, e.g. }$C_{2}C_{5}\not=0$\emph{\ if} $\nu ^{\prime }=n^{\prime
}=2<n.$

\bigskip

2.27. \textbf{Lemma}. \emph{In general we have}%
\begin{equation*}
n\nu \leq \mu _{\min }+n^{\prime }\nu ^{\prime }.\emph{\ }
\end{equation*}

\emph{If there is equality }$n\nu =\mu _{\min }+n^{\prime }\nu ^{\prime }$%
\emph{\ then either:}

\ \ \emph{(i) }$m=ln$\emph{\ (here }$\nu ^{\prime }=\nu -l(n-1)),$\emph{\ or}

\ \ \emph{(ii) }$m=ln+n^{\prime }$\emph{\ and }$\nu ^{\prime }=0.$

\emph{In the second case we have} 
\begin{equation}
\mu =\mu _{\min }\leq m(n-1)-n/2.  \tag{2.28}
\end{equation}

\bigskip

The numbers of double points hidden in places at infinity are estimated in
the next proposition, whose proof repeats the proof of Propositions 2.12 and
2.16 from [BZI]. Recall the notations $p^{\prime }=\gcd (p,q),$\emph{\ }$%
r^{\prime }=\gcd (r,s)$ and recall that $\delta _{\infty }=0$ if $p^{\prime
}=1$ and $\delta _{0}=0$ if $r^{\prime }=0.$

\bigskip

2.29. \textbf{Proposition.} (\emph{a) If }$ps\neq rq,$\emph{\ then}%
\begin{equation}
2\delta _{\infty }\leq p^{\prime }\nu _{\infty }\text{ if }p^{\prime
}>1,\;\;\;\;2\delta _{0}\leq r^{\prime }\nu _{0}\text{ if }r^{\prime }>1, 
\tag{2.30}
\end{equation}%
\emph{where }$\nu _{\infty }$\emph{\ (respectively }$\nu _{0})$\emph{\ is
the codimension of a corresponding stratum in the space of curves with
asymptotic (2.7) (respectively (2.8)).}

\emph{(b) If }$ps=rq,$\emph{\ then }%
\begin{equation}
2\delta _{\inf }\leq (p^{\prime }+r^{\prime })(\nu _{\inf }+1),  \tag{2.31}
\end{equation}%
\emph{where }$\nu _{\inf }=\nu _{0}+\nu _{\infty }+\nu _{\tan },$\emph{\ }$%
\nu _{0}$\emph{\ and }$\nu _{\infty }$\emph{\ are defined as in the point
(a), and }$\nu _{\tan }$\emph{\ is the number of first coinciding
coefficients of the Puiseux expansions of the branches }$\mathcal{C}_{\infty
}$\emph{\ and }$\mathcal{C}_{0}$\emph{\ which are not vanishing essential
Puiseux coefficients (}$\nu _{\tan }=u-$\emph{number of vanishing essential }%
$G_{j}$\emph{'s in (2.17)).}

\bigskip

This proposition admits the following improvements.

\bigskip

2.32. \textbf{Lemma}. \emph{(a) If }$ps=rq$\emph{\ and }$p^{\prime
}=r^{\prime }=1,$\emph{\ then }$2\delta _{\inf }\leq 2\nu _{\inf }$\emph{\
(here }$\nu _{0}=\nu _{\infty }=0).$

\emph{(b) If }$ps=rq$\emph{\ and }$\nu _{\inf }=0$\emph{\ or }$\nu _{\inf
}=1,$\emph{\ then }$2\delta _{\inf }\leq (p^{\prime }+r^{\prime })\nu _{\inf
}.$\emph{\ }

\emph{(c) Let }$ps\neq rq$\emph{\ and }$y=x^{q_{1}/p_{1}}(1+C_{l}t^{-l}+%
\ldots ),$\emph{\ }$C_{l}\neq 0,$\emph{\ as }$t\rightarrow \infty .$\emph{\
Then }$2\delta _{\infty }\leq 2\delta _{\infty ,\min }+p^{\prime \prime }\nu
_{\infty }^{\prime \prime },$\emph{\ where }%
\begin{equation*}
2\delta _{\infty ,\min }=(l-1)(p^{\prime }-1)+(p^{\prime \prime }-1),\emph{\ 
}
\end{equation*}%
$p^{\prime \prime }=\gcd (p,l)$\emph{\ and }$\nu _{\infty }^{\prime \prime }$%
\emph{\ is the corresponding codimension. (Analogous statement holds for }$%
2\delta _{0}).$

\bigskip

2.33. \textbf{Remark}. If $ps=rq,$ then%
\begin{equation*}
p=p_{1}p^{\prime },\;\;r=p_{1}r^{\prime },\;\;q=q_{1}p^{\prime
},\;\;s=q_{1}r^{\prime }.
\end{equation*}%
Here the property $\nu _{\tan }\geq 1$ means that 
\begin{equation*}
E_{1}^{p_{1}}=F_{1}^{p_{1}}
\end{equation*}%
in the Puiseux series $y=E_{1}x^{q_{1}/p_{1}}+\ldots $ and $%
y=F_{1}x^{q_{1}/p_{1}}+\ldots .$

\bigskip

\textbf{2.III. Spaces of parametric annuli.} We consider curves of the form $%
x=\varphi (t),$ $y=\psi (t)$, where 
\begin{equation}
\begin{array}{c}
\varphi =t^{p}+a_{1}t^{p-1}+\ldots +a_{p+r}t^{-r}, \\ 
\psi =t^{q}+b_{1}t^{q-1}+\ldots +b_{q+s}t^{-s}.%
\end{array}
\tag{2.34}
\end{equation}%
For fixed $p,r,q,s$ we denote by $\overline{Curv}=\overline{Curv}%
_{r,p;s,q}\simeq \mathbb{C}^{p+q+r+s}\setminus \left\{
a_{p+r}b_{q+s}\not=0\right\} $ the space of such curves.

It is easy to see that, upon application of a Cremona transformation and of
eventual change $t\rightarrow 1/t,$ we can divide all curves into the
following four types (recall that $p^{\prime }=\gcd (p,q),$ $r^{\prime
}=\gcd (r,s))$:%
\begin{equation}
\begin{array}{rl}
\text{Type }\binom{+}{+}: & \text{when }0<p<q,\text{ }0<r<s,\text{ }%
r^{\prime }\leq p^{\prime }\text{ and }\min \left( \frac{q}{p},\frac{s}{r}%
\right) \notin \mathbb{Z}; \\ 
\text{Type }\binom{-+}{+-}: & \text{when }0<q<p,\text{ }0<r<s\text{ and }%
p+r\leq q+s; \\ 
\text{Type }\binom{-}{+}: & \text{when }0<-r\leq p,\text{ }q>0,\text{ }s>0%
\text{ and }\frac{q}{p}\notin \mathbb{Z}; \\ 
\text{Type }\binom{-}{-}: & \text{when }0<-r\leq p,\text{ }0<-q\leq s\text{
and }p-\left\vert r\right\vert \leq s-\left\vert q\right\vert .%
\end{array}
\tag{2.35}
\end{equation}%
Graphically they are presented at the below figures. \medskip

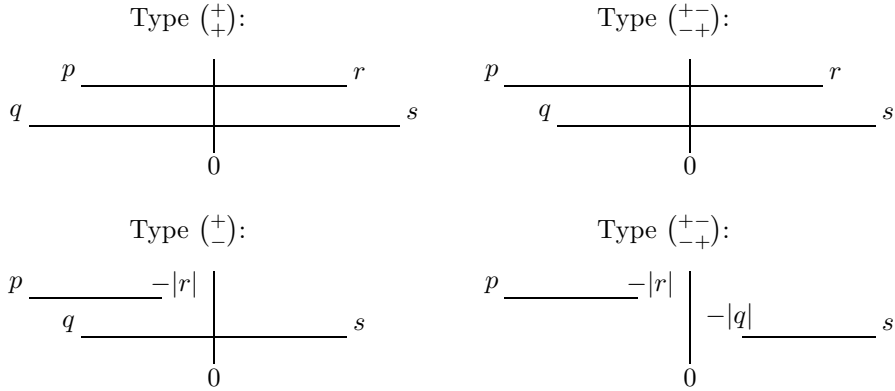
\begin{picture}(400,180)(0,20)
\put(5,130){\line(1,0){140}}
\put(25,145){\line(1,0){100}}
\put(75,120){\line(0,1){35}}
\put(0,135){\makebox(0,0){$q$}}
\put(150,135){\makebox(0,0){$s$}}
\put(20,150){\makebox(0,0){$p$}}
\put(130,150){\makebox(0,0){$r$}}
\put(75,115){\makebox(0,0){$0$}}
\put(65,170){\makebox(0,0){Type $\binom{+}{+}$:}}
\put(205,130){\line(1,0){120}}
\put(185,145){\line(1,0){120}}
\put(255,120){\line(0,1){35}}
\put(200,135){\makebox(0,0){$q$}}
\put(330,135){\makebox(0,0){$s$}}
\put(180,150){\makebox(0,0){$p$}}
\put(310,150){\makebox(0,0){$r$}}
\put(255,115){\makebox(0,0){$0$}}
\put(245,170){\makebox(0,0){Type $\binom{+-}{-+}$:}}
\put(25,50){\line(1,0){100}}
\put(5,65){\line(1,0){50}}
\put(75,40){\line(0,1){35}}
\put(20,55){\makebox(0,0){$q$}}
\put(130,55){\makebox(0,0){$s$}}
\put(0,70){\makebox(0,0){$p$}}
\put(60,70){\makebox(0,0){$-|r|$}}
\put(75,35){\makebox(0,0){$0$}}
\put(65,90){\makebox(0,0){Type $\binom{+}{-}$:}}
\put(185,65){\line(1,0){50}}
\put(275,50){\line(1,0){50}}
\put(255,40){\line(0,1){35}}
\put(270,58){\makebox(0,0){$-|q|$}}
\put(330,55){\makebox(0,0){$s$}}
\put(180,70){\makebox(0,0){$p$}}
\put(240,70){\makebox(0,0){$-|r|$}}
\put(255,35){\makebox(0,0){$0$}}
\put(245,90){\makebox(0,0){Type $\binom{+-}{-+}$:}}
\end{picture}

The space $\overline{Curv}$ (for fixed $p,r,q,s)$ admits action of a group $%
\mathcal{G}$ generated by:

\begin{itemize}
\item the multiplication of $t$ by $\lambda ^{-1}$ accompanied with
multiplication of $x$ by $\lambda ^{p}$ and of $y$ by $\lambda ^{q};$

\item the addition of a constant to $x$ (respectively to $y)$ if $r>0$
(respectively if $q>0);$

\item the change $y\rightarrow y+P(x)$ for a polynomial $P$ of degree
\end{itemize}

\begin{equation}
k=\min \left( \left[ \frac{q}{p}\right] ,\left[ \frac{s}{r}\right] \right) 
\text{ for }\binom{+}{+},\;=\left[ \frac{q}{p}\right] \text{ for }\binom{-}{+%
},\text{ }=0\text{ for }\binom{-+}{+-}\text{ and }\binom{-}{-}.  \tag{2.36}
\end{equation}

\medskip

2.37. \textbf{Definition.} The space $Curv=\overline{Curv}/\mathcal{G}$ is
called the \textit{space of annuli}. Its is a quasi-projective variety of
dimension%
\begin{equation}
\sigma :=\dim Curv=p+q+r+s-1-\varepsilon -k,  \tag{2.38}
\end{equation}%
where $\varepsilon =2$ for Type $\binom{+}{+}$ and Type $\binom{-+}{+-},$ $%
\varepsilon =1$ for Type $\binom{-}{+}$ and $\varepsilon =0$ for Type $%
\binom{-}{-}.$

\bigskip

A typical element from $Curv$ has $\delta _{\max }$ simple double points. If 
$\xi \in Curv$ is such that its image $\mathcal{C}=\xi (\mathbb{C}^{\ast })$
does not have self-intersections, then its double points are hidden either
at singular points $\xi (t_{1}),\ldots ,\xi (t_{N})$ or at infinity.

\bigskip

2.39. \textbf{Definition.} A point $t_{j}$ is \textit{singular} for a
parametric curve $\xi :\mathbb{C}^{\ast }\rightarrow \mathbb{C}^{2}$ iff $%
\xi ^{\prime }(t_{j})=0;$ (therefore a self-intersection of smooth branches
of $\mathcal{C}=\xi (\mathbb{C}^{\ast })$ is not regarded as singular point
of the immersion $\xi ,$ though it is a singular point of $\mathcal{C}$). We
have 
\begin{equation*}
\varphi (t)=x_{j}+(t-t_{j})^{n_{j}}(\kappa _{j}+\ldots ),\;\;\psi
(t)=y_{j}+O((t-t_{j})^{2}),
\end{equation*}%
where $n_{j}$ is called the $x$-\textit{order of} $t_{j}.$ The singular
point $t_{j}$ is characterized by its $y$-\textit{codimension} $\nu _{j}$
(in the sense of Proposition 2.20) and by its Milnor number $\mu
_{t_{j}}=2\delta _{t_{j}}.$ We define the \textit{external codimension of} $%
t_{j}$ as 
\begin{equation*}
ext\nu _{j}=(n_{j}-2)+\nu _{j}.
\end{equation*}%
(Note that $n_{j}-2$ is the number of conditions that $\varphi ^{\prime }(t)$
has (somewhere) zero of order $n_{j}-1.)$

The above notions of $x$-order and of $y$-codimension are not symmetric with
respect to the change $x\longleftrightarrow y.$ In fact, we use these
notions (in this form) when $p+r\leq q+s;$ so Types $\binom{+}{+},$ $\binom{%
-+}{+-}$ and $\binom{-}{-}$ are included here. But for Type $\binom{-}{+}$
with $q+s<p-\left\vert r\right\vert $ we define $n_{j}$ as the $y$-order and 
$\nu _{j}$ as the $x$-codimension of the singularity.

When some double points are hidden at infinity, then the corresponding
external codimensions are 
\begin{equation*}
ext\nu _{0}=\nu _{0},\;\;ext\nu _{\infty }=\nu _{\infty },\;\;ext\nu _{\inf
}=\nu _{\inf },
\end{equation*}%
where $\nu _{0},$ $\nu _{\infty }$ and $\nu _{\inf }=\nu _{0}+\nu _{\infty
}+\nu _{\tan }$ were defined in Proposition 2.29 (with the agreement $\nu
_{\tan }=0$ when $ps\not=rq).$

The space $\overline{Curv}$ contains curves such that the map $\xi :\mathbb{C%
}^{\ast }\rightarrow \mathcal{C}$ is several-to-one. Such curves are called 
\textit{multiply covered} (or \textit{non-primitive}) and form an algebraic
subvariety $Mult$ of $\overline{Curv}.$ If $\xi \in Mult$ then some its
singular points have infinite codimension.

\medskip

The conjecture following was mentioned in Introduction. It is crucial in the
sequel sections and is proved in the sequent paper [BZIII].

\bigskip

2.40. \textbf{Conjecture (codimension bound)}. \emph{Suppose that }%
\begin{equation}
ext\nu _{\inf }+\sum_{j=1}^{N}ext\nu _{t_{j}}\leq \sigma =\dim Curv. 
\tag{2.41}
\end{equation}%
\emph{Then the degenerations as described in Definition 2.39 occur along an
algebraic subvariety of the space }$\overline{Curv}\setminus Mult$\emph{\ of
codimension }$ext\nu _{\inf }+\sum ext\nu _{t_{j}}.$\emph{\ If the above sum
of external codimensions is greater than }$\sigma ,$\emph{\ then it equals }$%
\infty $\emph{\ and the degeneration occurs along a subvariety consisting of
multiply covered curves.}

\emph{The inequality (2.41) is called the \textbf{regularity condition}.}

\bigskip

2.42. \textbf{Remark.} The \textit{multiply covered curves} (or
non-primitive curves) are such that the map $\xi :\mathbb{C}^{\ast
}\rightarrow \mathbb{C}^{2}$ is several-to-one. By the L\"{u}roth theorem
(or the Stein factorization) any multiply covered curve has the form $\xi =%
\tilde{\xi}\circ \omega ,$ where either $\omega (t)=t^{d}$ and $\tilde{\xi}$
is primitive, or $\omega :\mathbb{C}^{\ast }\rightarrow \mathbb{C}$ is a
Laurent polynomial and the mapping $\tilde{\xi}:\mathbb{C}\rightarrow 
\mathbb{C}^{2}$ is polynomial and primitive.

\bigskip

2.43. \textbf{Remark.} In [BZIII] we introduce (following Orevkov) a
so-called rough M-number $\overline{M}_{z}$ of a cuspidal singular point $P$
of the curve $\mathcal{C}.$ When the curve has the form (2.4) and $n$ is the
multiplicity of $P=(0,0)$, i.e. the degree of the first nonzero term of the
Taylor expansion at $P$ of the polynomial defining $\mathcal{C}$, then $%
\overline{M}_{P}$ equals the external codimension of the singularity $ext\nu
_{P}.$ Otherwise $\overline{M}_{P}<ext\nu _{P},$ but the difference is well
controlled.

We prove in [BZIII] the following bounds:%
\begin{equation*}
ext\nu _{\inf }+\sum \overline{M}_{P_{i}}\leq p+q+r+s+1-K,\text{ \ \ }K=\min
\left( [q/p],\left[ s/r\right] \right)
\end{equation*}%
for Types $\binom{+}{+}$ and $\binom{-+}{+-},$ 
\begin{equation*}
ext\nu _{\inf }+\sum \overline{M}_{P_{i}}\leq p+q-|r|+s+2-K+[(|r|-1)/s]
\end{equation*}%
for Type $\binom{-}{+}$ and 
\begin{equation*}
ext\nu _{\inf }+\sum \overline{M}_{P_{i}}\leq
p-|q|-|r|+s+3+[(|r|-1)/s]+[(|q|-1)/p]
\end{equation*}%
for Type $\binom{-}{-}$.

\bigskip

\textbf{2.IV. Handsomeness.} The division of annuli into the four types in
(2.35) is not completely precise. It depends on the choice of the reducing
automorphism of $\mathbb{C}^{2}.$ For example, if $x=t^{2}+\ldots +t^{-3},$ $%
y=t^{6}+\ldots +t^{-4}$ is of Type $\binom{+}{+}$ then the change $%
y\rightarrow y-x^{3}=t^{\tilde{q}}+\ldots +t^{-9}$ may give a curve of Type $%
\binom{+}{+}$ or of Type $\binom{-+}{+-}$ or of Type $\binom{-}{+}.$ In
order to avoid this ambiguity we introduce the notion of handsome curve,
which also will turn out useful in estimates in the further sections.

\bigskip

2.43. \textbf{Definition.} A curve $\xi $ (of one of the four types in
(2.35)) is called \textit{non-handsome} if either

\begin{itemize}
\item $\frac{q}{p}\in \mathbb{Z}$ and $r<p$ for Type $\binom{+}{+},$ or

\item $\frac{p}{q}\in \mathbb{Z}$ and $s<q$ or $\frac{s}{r}\in \mathbb{Z}$
and $p<r$ for Type $\binom{-+}{+-},$ or

\item $\frac{p}{q}\in \mathbb{Z}$ and $s<q$ for Type $\binom{-}{+}.$
\end{itemize}

Otherwise the curve is \textit{handsome}.

\bigskip

2.44. \textbf{Proposition.} \emph{Any non-handsome curve can be transformed
using a Cremona automorphism and/or the change }$t\rightarrow 1/t$\emph{\ to
a handsome curve (of one of the types in (2.35)).}\medskip

\textit{Proof}. 1. Suppose firstly that a curve of Type $\binom{-+}{+-}$ is
non-handsome. Assume that $\frac{p}{q}$ is integer and $s<q=p^{\prime };$
the case of $\frac{s}{r}$ integer and $p<r=r^{\prime }$ is treated
analogously. We have $x=t^{p_{1}q}+\ldots +at^{-r},$ $y=t^{q}+\ldots
+bt^{-s}.$ We apply the changes $x\rightarrow x-$const$\cdot y^{l}$ as many
times as possible (in order to diminish $\deg x)$. We obtain $x=t^{\hat{p}%
}+\ldots +ct^{-p_{1}s},$ $y=t^{q}+\ldots +bt^{-s},$ where either (i) $1<\hat{%
p}/q\notin \mathbb{Z},$ or (ii) $\hat{p}<0,$ or (iii) $0<\hat{p}/q<1$.

In the case (i) we have Type $\binom{+}{+}$ (after swapping $x$ with $y$)
with the exponents $\tilde{p}=q,$ $\tilde{r}=s,$ $\tilde{q}=\hat{p},$ $%
\tilde{s}=p_{1}s.$ It would be non-handsome only when $\tilde{r}>\tilde{p},$
i.e. $s>q$ (here $t\rightarrow 1/t);$ but we have assumed reverse inequality.

In the case (ii) we get a curve of Type $\binom{-}{+}$ with $\tilde{p}%
=p_{1}s,$ $\tilde{r}=\hat{p},$ $\tilde{q}=s,$ $\tilde{s}=q.$ Here $\tilde{p}%
^{\prime }=s<\tilde{s}=p^{\prime }$ (by assumption) and the curve is
handsome.

In the case (iii) we get a curve of Type $\binom{-+}{+-},$ with exponents
like in the case (i). It would be non-handsome if either $r>\tilde{p}$ (it
is not the case), or $\tilde{s}/\tilde{r}=q/\hat{p}\in \mathbb{Z}$ and $%
\tilde{r}=\hat{p}>\tilde{q}=s.$ This is the same situation we started with,
but now $\tilde{q}<q$ and $\tilde{s}>s.$ Repeating the above reduction
process several times we must arrive to one of the cases (i), (ii) or (iii)
with $q\leq s$ (handsomeness).\medskip

2. Suppose that we have a non-handsome curve of Type $\binom{-}{+},$ i.e. $%
\frac{p}{q}$ is integer and $s<q$ (of course, also $r<0).$ Then after
application of transformations like in the point 1 we get a curve like in
the point 1 and with the same three possibilities. The further proof is also
the same.\medskip

3. Suppose that we have a non-handsome curve of Type $\binom{+}{+}.$ Thus $%
\frac{q}{p}$ is integer and $p=p^{\prime }>r:$ $x=t^{p}+\ldots +at^{-r},$ $%
y=t^{q_{1}p}+\ldots +bt^{-s}$ (where $\frac{s}{r}<q_{1}).$ The changes $%
y\rightarrow y+$const$\cdot x^{l}$ reduce it to the form $x=t^{p}+\ldots
+at^{-r},$ $y=t^{\hat{q}}+\ldots +ct^{-q_{1}r}.$ We have three
possibilities: (i) $1<\hat{q}/p\notin \mathbb{Z},$ (ii) $\hat{q}<0,$ (iii) $%
0<\hat{q}/p<1.$

In the case (i) we get a curve of Type $\binom{+}{+}$ which would be
non-handsome iff $r>p$ (it is not the case).

In the case (ii) we get a curve of Type $\binom{-}{+}$ with $\tilde{p}%
=q_{1}r,$ $\tilde{r}=\hat{q}<0,$ $\tilde{q}=r,$ $\tilde{s}=p.$ It would be
non-handsome iff $\tilde{q}=r>\tilde{s}=p,$ but it is not the case.

In the case (iii) we get a curve of Type $\binom{-+}{+-}$ and further proof
runs along the lines of the point 1. $\qed$

\bigskip

2.46. \textbf{Remark.} Assuming that a curve is handsome we can try to apply
changes like in the proof of the latter proposition. It turns out that
either the transformed curve falls out of the list in (2.35) or becomes
non-handsome. To be correct, this statement does not apply to the cases $%
\frac{q}{p}=\frac{q}{r}\in \mathbb{Z}$ and $\frac{s}{r}=\frac{s}{p}\in 
\mathbb{Z}.$ In this sense the list of the four types $\binom{+}{+}$, $%
\binom{-+}{+-},$ $\binom{-}{+}$ and $\binom{-}{-}$ of handsome curves is
complete and unique.

On the other hand, the handsome curves are such that the codimension of
their degenerations at infinity are smallest possible.

It is a good place to say why all the curves from the list in Main Theorem
are pairwise different. The division into the four types of handsome curves
is a preliminary classification; essentially it is a classification with
respect to the leading terms of the expansions of the curves as $%
t\rightarrow \infty $ and as $t\rightarrow 0.$ The further classification
within a fixed type follows from different types of finite singularities
and/or from different details of the expansions at $t=\infty $ and at $t=0.$

In the below proofs we concentrate on detection of the cases of embedding of 
$\mathbb{C}^{\ast }$ omitting the details of the analysis which cases are
really different.

\bigskip

\textbf{2.V. Scheme of the proof}. We introduce the quantity%
\begin{equation}
\mathcal{E}=n_{\infty }\nu _{\infty }+n_{0}\nu _{0}+\sum_{j=1}^{N}n_{j}\nu
_{j},\;\;\;\;\mathrm{if}\;\;\;\;ps\neq rq,  \tag{2.47}
\end{equation}%
or%
\begin{equation}
\mathcal{E}=n_{\inf }(\nu _{\inf }+1)+\sum_{j=1}^{N}n_{j}\nu _{j},\;\;\;\;%
\mathrm{if}\;\;\;\;ps=rq.  \tag{2.48}
\end{equation}%
It should satisfy the inequality 
\begin{equation}
2\delta _{\max }\leq \mathcal{E}  \tag{2.49}
\end{equation}%
(by (2.16), (2.21), (2.30) and (2.31)). We shall strive to detect the cases
when the inequality (2.48) holds true. We shall use the \textit{reserve} 
\begin{equation*}
\Delta :=2\delta _{\max }-\mathcal{E}.
\end{equation*}%
If $\Delta >0,$ then the curve $\xi $ is not an embedding.

If $\Delta =0$, and cannot be negative, then we say that the case is \textit{%
strict}; it means that all inequalities leading to $\Delta \leq 0$ must be
equalities. For instance, $\mu =n\nu $ for all finite singular points.

In estimation of $\mathcal{E}$ from the above we use some natural
restrictions. For example, in Type $\binom{+}{+}$ with $ps\neq rq$ we have 
\begin{equation}
\begin{array}{c}
\sum_{j=1}^{N}(n_{j}-1)\leq p+r,\;\;n_{j}\leq p+r, \\ 
\nu _{\infty }+\nu _{0}+\sum_{j=1}^{N}\nu _{j}+\sum (n_{j}-2)\leq \sigma .%
\end{array}
\tag{2.50}
\end{equation}%
Like in [BZI] one shows that the maximum of $\mathcal{E}$ is achieved in the
case when only one singular point is not a standard cusp. (The analysis is
slightly different for Type $\binom{-}{-},$ where the double points hide
rather at infinity.)

Next, after detecting some genuine cases of $\mathbb{C}^{\ast }$-embeddings,
one arrives to the situation with only one singular point, which can be put
at $t_{1}=1.$ Moreover, the $x$-order $n=n_{1}$ of this point can equal $%
p+r-1$ or $p+r.$ Here, in order to get more precise estimate of $2\delta
_{1} $ one has to consider curves of the form 
\begin{equation}
\varphi =(t-1)^{n}P(t)t^{-r},\;\;\psi =(t-1)^{m}Q(t)t^{-s}  \tag{2.51}
\end{equation}%
(e.g. of Type $\binom{+}{+})$. Here the bound (2.23) is used, $2\delta
_{1}\leq \mu _{\min }+n^{\prime }\nu ^{\prime }.$ But the sum of
codimensions $\nu ^{\prime }+\nu _{0}+\nu _{\infty }$ cannot be estimated
directly by the dimension of the space of curves of the form (2.51), i.e. by 
$\deg P+\deg Q-\tilde{k},$ where $\tilde{k}$ counts the changes $%
y\rightarrow y+$\textrm{const}$\cdot x^{l}$ preserving (2.51).

\bigskip

2.52. \textbf{Proposition.} \emph{Assume that a non-primitive curve, which
satisfies the regularity condition (2.41), has the form (2.51) with singular
points }$t_{1}=1,t_{2},\ldots ,t_{N}$\emph{. Then we have }%
\begin{equation*}
\nu _{\inf }+\nu ^{\prime }+ext\nu _{2}+\ldots +ext\nu _{N}\leq \deg P+\deg
Q-k+[(m-1)/n],
\end{equation*}%
\emph{where }$k$\emph{\ is the same as in (2.36) and }$[\cdot ]$\emph{\
denotes the integer part. }\medskip

\textit{Proof}. (This proof repeats the proof of Lemma 3.9 from [BZI].) We
consider curves of the type (2.51).They are defined by vanishing of the
first $m-1$ Puiseux coefficients $C_{1}^{(1)},\ldots ,C_{m-1}^{(1)},$ where $%
[\frac{m-1}{n}]$ of them are essential (with the lower indices being
multiples of $n).$ Therefore 
\begin{equation*}
\nu _{1}=[(m-1)/n]+\nu _{1}^{\prime },
\end{equation*}%
where $\nu _{1}^{\prime }$ is the codimension in the space of curves (2.51).

Now the proposition follows from the inequality (2.41). $\qed$

\bigskip

\section{Annuli of Type $\binom{-}{-}$}

We begin the proof of Main Theorem with the type $\binom{-}{-},$ because
many situations with curves of other types are reduced to Type $\binom{-}{-}%
. $

Recall that we deal with curves of the form 
\begin{equation*}
\varphi =t^{p}+\ldots +a_{p-\left\vert r\right\vert }t^{\left\vert
r\right\vert },\;\;\;\;\psi =t^{-\left\vert q\right\vert }+\ldots
+b_{s-\left\vert q\right\vert }t^{-s},
\end{equation*}%
where we additionally assume that (see (2.35)) 
\begin{equation}
p-\left\vert r\right\vert \leq s-\left\vert q\right\vert .  \tag{3.1}
\end{equation}%
The further analysis is divided into four cases:%
\begin{equation*}
p-\left\vert r\right\vert =0,\;\;p-\left\vert r\right\vert
=1,\;\;p-\left\vert r\right\vert =2,\;\;p-\left\vert r\right\vert \geq 3.
\end{equation*}

\bigskip

\textbf{3.I. The case }$p=\left\vert r\right\vert .$ Therefore $\varphi
=t^{p},$ $p>0.$ The following lemma can be also found in [Ka] (with a
slightly different formulation).

\bigskip

3.2. \textbf{Lemma.} \emph{Any annulus of the form }$x=t^{p},$\emph{\ }$%
y=\psi (t),$\emph{\ where }$\psi =b_{0}t^{q}+b_{1}t^{q-1}+\ldots
+b_{q+s}t^{-s}$\emph{, can be reduced to }$x=t^{p},$\emph{\ }$y=t^{d}+\gamma
_{1}t^{-p}+\ldots +\gamma _{l}t^{-lp},$\emph{\ where }$\gcd (p,d)=1.$\emph{\
It is item (a) of Main Theorem, where }$d$\emph{\ can be either positive or
negative. }\medskip

\textit{Proof}. If $p=1$ then $\psi $ can be reduced to a polynomial of $%
1/t. $

Let $p>1.$ The double point equations $\varphi (t^{\prime })=\varphi (t),$ $%
\psi (t^{\prime })=\psi (t)$ cannot have solutions $t^{\prime },t\in \mathbb{%
C}^{\ast }.$ Since $t^{\prime }=\zeta t,$ $\zeta \neq 1$ a root of unity of
degree $p,$ we get the equation 
\begin{equation*}
b_{0}(\zeta ^{q}-1)t^{q}+b_{1}(\zeta ^{q-1}-1)t^{q-1}+\ldots +b_{q+s}(\zeta
^{-s}-1)t^{-s}=0.
\end{equation*}%
For each $\zeta $ the monomial in the left-hand side can contain at most one
monomial. If all these monomials vanished, then $\psi $ would depend on $%
t^{p}$ and the curve would be multiply covered. Therefore only one monomial $%
b_{d}t^{d}$ is such that $\zeta ^{d}\neq 1$ for all $\zeta ,$ all other
monomials in $\psi $ are powers of $t^{p}.$ The positive powers of $t^{p}$
can be killed, but the negative powers of $t^{p}$ remain. $\qed$

\bigskip

\textbf{3.II. The case }$p=\left\vert r\right\vert +1.$ By rescaling $t$ we
can assume 
\begin{equation}
\varphi =(t-1)t^{\left\vert r\right\vert },\;\;\;\;\psi =Q(1/t)  \tag{3.3}
\end{equation}%
for a polynomial $Q$ of degree $s$ and $\mathrm{ord}_{t=0}Q=\left\vert
q\right\vert \leq s-1.$

The following constructions are important.

\bigskip

3.4. \textbf{Definition.} Suppose that we have a $\mathbb{C}^{\ast }$%
-embedding $\xi =(\varphi ,\psi )$ such that 
\begin{equation}
\varphi =(t-t_{1})^{n}t^{-r},\;\;\;\;n=p+r.  \tag{3.5}
\end{equation}%
Then the curve 
\begin{equation}
\tilde{\xi}=(\tilde{\varphi},\tilde{\psi}):=(\varphi ,\varphi \psi ), 
\tag{3.6}
\end{equation}%
is said to be obtained from $\xi $ by the \textit{tower transformation}.

We have also analogous tower transformation $(\varphi ,\psi )\rightarrow
(\varphi \psi ,\psi )$ when $\psi =$\linebreak $(t-t_{1})^{m}t^{-s},$ but we
shall mainly use the transformation (3.6).

In the case of curves of the form (3.3), i.e. with $n=1$ and $|r|>0,$ we
define the \textit{reverse tower transformation} 
\begin{equation}
\mathcal{T}:\xi \rightarrow \hat{\xi}=(\hat{\varphi},\hat{\psi}):=(\varphi ,%
\left[ \psi (t)-\psi (t_{1})\right] /\varphi (t)).  \tag{3.7}
\end{equation}

Of course, the tower transformation and the reverse tower transformation can
be regarded as the standard blowing down and blowing up constructions from
the projective algebraic geometry.\medskip

The following result does not require proof.

\bigskip

3.8. \textbf{Lemma.} \emph{(a) If }$\xi $\emph{\ is a }$\mathbb{C}^{\ast }$%
\emph{-embedding, then }$\tilde{\xi}$\emph{\ and }$\hat{\xi}$\emph{\ are
also }$\mathbb{C}^{\ast }$\emph{-embeddings. }

\emph{(b) }$\mathcal{T}$\emph{\ is the right inverse transformation to the
map }%
\begin{equation*}
\mathcal{T}^{-1}:\left( \varphi ,\psi \right) \rightarrow (\varphi ,\varphi
\psi +K)
\end{equation*}%
\emph{(where }$K$\emph{\ is a constant), which is equivalent to }$\xi
\rightarrow \tilde{\xi}$\emph{. If }$q>0,$\emph{\ then the constant }$K$%
\emph{\ is not defined uniquely; but when }$q<0$\emph{\ and }$\varphi \psi $%
\emph{\ does not have pole at }$t=\infty ,$\emph{\ then we put }%
\begin{equation*}
K=-\left( \varphi \psi \right) (\infty ).
\end{equation*}

\emph{(c) If }$n=2$\emph{\ and }$|r|=1,$\emph{\ then }$\varphi
+1/4=(t-1/2)^{2}$\emph{\ and we can apply another `tower transformation',
which takes the form }%
\begin{equation}
\left( \varphi ,\psi \right) \rightarrow \left( \varphi ,\left( \varphi
+1/4\right) \psi \right)  \tag{3.9}
\end{equation}%
\emph{and which is not equivalent to }$\tilde{\xi}.$

\bigskip

Let us return to annuli of the form (3.3). We shall apply to them the tower
transformations and the reverse tower transformations. However Lemma 3.8(c)
shows that the case with $\left\vert r\right\vert =1$ should be treated
separately.

\bigskip

3.10. \textbf{Lemma.} \emph{Any annulus of Type }$\binom{-}{-}$\emph{\ or }$%
\binom{-}{+}$\emph{\ with }$\varphi =t(t-1)$\emph{\ is equivalent to an
annulus obtained from }$\psi _{0}=(t-\frac{1}{2})$\emph{\ by applying:}

\begin{itemize}
\item $m$ \emph{times the operation (3.9), }$m=0,1,2,\ldots ,$ \emph{and}

\item $s$ \emph{times the reverse tower transformation (3.7), }$s=1,2,\ldots
.$
\end{itemize}

\emph{This gives the series (b) of Main Theorem. }\medskip

\textit{Proof}. The function $\psi (t)$ has pole at $t=0$ of order $s>0.$ We
apply the tower transformation (3.6) several times, just to reduce the pole
at $t=0:$ $\psi \rightarrow \psi _{1}=\psi \varphi ^{s}.$ The polynomial
curve $(\varphi ,\psi _{1})$ does not have double points. Like in Lemma 3.2
we find that $\psi _{1}=$const$\cdot (t-1/2)^{2m+1}+$(polynomial in $\varphi 
$) for some $m\geq 0.$ After a normalization we can assume that $\psi
_{1}=\psi _{2}+L(\varphi ),$ where $\psi _{2}=(t-1/2)^{2m+1}$ and $L$ is a
polynomial of degree $\leq s.$

We see that $\psi _{2}=(\varphi +1/4)^{m}\psi _{0}.$ Application of the
reverse transformation equation $\psi _{1}\rightarrow \psi =\psi
_{1}/\varphi ^{s}$ is the same as application of the reverse tower
transformation (3.7) to $\psi _{2}$. Hence we obtain the series (b) of Main
Theorem.

Note that application of $\mathcal{T}$ to $\xi $ with $\psi =t-1/2$ gives $%
\psi =1/t$ and to $\xi $ with $\psi =1/t$ gives $\psi =-1/t^{2};$ these two
cases are included into the series (a) of Main Theorem. Next, application of 
$\mathcal{T}$ to $\xi $ with $\psi =(t-1/2)^{3}$ gives $\psi =t+\frac{1}{4}%
t^{-1}+$const$,$ which belongs to the series (j) of Main Theorem. $\qed$

\bigskip

Assume now that%
\begin{equation*}
\left\vert r\right\vert \geq 2.
\end{equation*}%
By applying iterations of the reverse tower transformation to one annulus we
obtain a series of annuli. Therefore our task is to determine the initial
term of any such series.

By Lemma 3.8(b) application of $\mathcal{T}^{-1}$ gives a curve of Type $%
\binom{-}{-},$ and $\mathcal{T}^{-1}$ is uniquely defined here, only when $%
\left( \varphi \psi \right) (t)$ is a polynomial in $1/t,$ i.e. when $%
p-\left\vert q\right\vert \leq 0.$ Note, however, that the curve $\mathcal{T}%
^{-1}\xi $ may not satisfy the property (3.1). Anyway, we have to consider
the following two possibilities for the initial curve of such a series:%
\begin{equation*}
\text{A. }\left\vert q\right\vert =s,\;\;\;\text{B. }\left\vert q\right\vert
<p,\text{ }|q|<s.
\end{equation*}%
\medskip

\textbf{Case A.} Lemma 3.2 implies that either $p=ls$ or $\left\vert
r\right\vert =ls.$ These two possibilities lead to two series $\left(
\varphi ,\psi _{j}\right) $ defined by $\psi _{0}(t)=t^{-s},$ $\psi
_{j+1}(t)=\left[ \psi _{j}(t)-\psi _{j}(1)\right] /\varphi (t).$

\bigskip

3.11. \textbf{Lemma}. \emph{In this way we obtain items (c) and (d) of Main
Theorem}.

\bigskip

\textbf{Case B.} We need some estimates. We treat separately three
possibilities:%
\begin{equation*}
\text{B.1. }p^{\prime }\geq r^{\prime },2;\;\;\text{B.2. }r^{\prime }\geq
2,p^{\prime };\;\;\text{B.3. }p^{\prime }=r^{\prime }=1<n=2.
\end{equation*}%
Note that $n=2$ is the maximal $x$-order of eventual (unique) singular point
at $t_{1}=(p-1)/p$ and $p^{\prime }=\gcd (p,q),$ $r^{\prime }=\gcd (r,s).$

\bigskip

\textbf{B.1.} Accordingly to Subsection 2.III we have $2\delta _{\infty
}+2\delta _{0}+2\delta _{t_{1}}\leq \mathcal{E}=n_{\infty }\nu _{\infty },$
where $n_{\infty }=p^{\prime }$ and $\nu _{\infty }=\sigma =\dim
Curv=s-\left\vert q\right\vert $ (see (2.41)). Here $2\delta _{\infty
}+2\delta _{0}+2\delta _{t_{1}}$ should equal $2\delta _{\max }=0\cdot
(s-\left\vert q\right\vert -1)-(p^{\prime }+r^{\prime
}-1)+(ps-|r||q|)=p(s-\left\vert q\right\vert )+\left\vert q\right\vert
-p^{\prime }-r^{\prime }+1$ (see (2.11)). Hence the reserve%
\begin{equation}
\Delta =2\delta _{\max }-\mathcal{E}  \tag{3.12}
\end{equation}%
should be non-positive (compare Subsection 2.IV). But we have 
\begin{equation*}
\Delta =(p-p^{\prime })(s-\left\vert q\right\vert )-r^{\prime }+(\left\vert
q\right\vert -p^{\prime })+1.
\end{equation*}%
Since $p^{\prime }=\gcd (p,\left\vert q\right\vert )$ and $\left\vert
q\right\vert <p,$ we have $p-p^{\prime }\geq p^{\prime }.$ Since also $%
p^{\prime }\geq r^{\prime },$ we find $\Delta \geq (p^{\prime }-r^{\prime
})+(\left\vert q\right\vert -p^{\prime })+1>0.$ So this case is not realized.

\bigskip

\textbf{B.2.} We perform analogous calculations as in the point B.1 and we
get $\mathcal{E}=n_{0}\nu _{0}=r^{\prime }(s-\left\vert q\right\vert )$ and 
\begin{equation*}
\Delta =\left( |r|-r^{\prime }\right) \left( s-\left\vert q\right\vert
\right) +\left( s-\left\vert r^{\prime }\right\vert \right) +1-p^{\prime }.
\end{equation*}%
In order that $\Delta \leq 0$ it should be%
\begin{equation*}
\left\vert r\right\vert =r^{\prime }=s,
\end{equation*}%
but the case is not strict (i.e. we can get $\Delta <0,$ see Subsection
2.IV). We have $\varphi =t^{s}(-1+t),$ $\psi =t^{-s}(b_{0}+b_{1}t+\ldots
+b_{s-\left\vert q\right\vert }t^{s-\left\vert q\right\vert }).$ Therefore
the curve $t\rightarrow \left( x,y\right) =\left( \varphi (t),\varphi
(t)\psi (t)\right) $ is a polynomial curve without self-intersections. After
the change $y\rightarrow y+b_{0}$ we get the curve%
\begin{equation}
x=t^{s}(t-1),\;\;\;y=t^{d}Q(t),\;\;\;d\geq 1,  \tag{3.13}
\end{equation}%
where the polynomial $Q$ is of degree $e$ and satisfies $Q(0)\neq 0.$
Moreover, $\deg y=d+e=s+1-|q|<\deg x.$

Here $d=\mathrm{ord}_{t=0}y$ should be $>1$ when $\nu _{0}>1$ (this must
occur when $r^{\prime }>n_{1}=2).$ But when $r^{\prime }=2$ the double
points may hide themselves at the unique singular point $t_{1}=s/(s+1).$ We
have then two possibilities:%
\begin{equation*}
(i)\;\;d>1;\;\;\;\;\;(ii)\;\;d=1.
\end{equation*}

\textit{Subcase} $\left( i\right) .$\medskip

3.14. \textbf{Lemma.} \emph{If }$d>1,$\emph{\ then }$Q(t)\equiv $const\emph{%
\ in (3.13). This implies that either }$\varphi =t^{ld}(t-1)$\emph{\ or }$%
\varphi =t^{ld-1}(t-1)$\emph{\ and the reverse tower transformations produce
the series (e) and (f) respectively of Main Theorem. }\medskip

\textit{Proof}. Since the curve (3.13) is contractible and singular, we can
use the Zaidenberg--Lin theorem [ZaLi]. It says that there exists a
composition of elementary transformations reducing (3.13) to a
quasi-homogeneous curve $(\alpha t^{u},\beta t^{v}),$ $\gcd (u,v)=1$.
Moreover, if $u=\deg \varphi $ and $v=\deg \psi $ are relatively prime then
the reduction to the quasi-homogeneous curve uses only only translations and
elementary transformations which do not change the degrees of $\varphi $ and 
$\psi .$

It implies that for $\gcd (p,d+e)=1,$ $p=s+1,$ the curve is
quasi-homogeneous. Therefore we must consider the case 
\begin{equation*}
p=l(d+e).
\end{equation*}%
Moreover, the (unique) singular point $t=0$ has only one characteristic
pair. Therefore either ($\alpha $) $\gcd (s,d)=1,$ or ($\beta $) $r=md.$

In the case ($\alpha $) the characteristic pair of the singularity is $%
(s,d). $ By the Zai\-den\-berg--Lin theorem the same is the characteristic pair
at $t=\infty ,$ after applying elementary transformations resulting in
relatively prime degrees of the components. But the first elementary
transformation should be $x\rightarrow x_{1}=x-M(y)$ for a polynomial $M$ of
degree $l$ (and without the constant term). This is not the only change. The
next one should be of the type $y\rightarrow y+N(x);$ thus $\deg x_{1}\leq
\deg y<s.$ It follows that the characteristic pair $\left( u,v\right) ,$ $%
u,v>1,$ at $t=\infty $ should satisfy $\max (u,v)<\max (s,d)=s.$ So this
case does not occur.

Consider the case ($\beta $). We have $2\delta _{\max }=\left( s-1\right)
(d+e-1)$ (see (2.10) and [BZI]), which should equal $2\delta _{0}+2\delta
_{\infty }.$ Here $2\delta _{0}=\mu _{0}\leq \mu _{0,\min }+d\nu _{0},$ $\mu
_{0,\min }=s(d-1)$ and $2\delta _{\infty }\leq (d+e)\nu _{\infty },$ where $%
\nu _{0},\nu _{\infty }$ are the codimensions at $t=0,\infty $ (see Section
2.II). The codimensions $\nu _{0,\infty }$ should be considered as the
codimensions in the $e$-dimensional space of curves (3.13). Here one applies
an analogue of Proposition 2.52, which states that $\nu _{0}+\nu _{\infty
}\leq \deg Q-[\frac{\deg x}{\deg y}]=e-l+m.$ Of course, the maximum of $%
\mathcal{E}=\mu _{0,\min }+d\nu _{0}+(d+e)\nu _{\infty }$ is achieved when
one puts $\nu _{0}=0$ and $\nu _{\infty }=e-l+m;$ then $\mathcal{E}%
=s(d-1)+(d+e)\left( e-l+m\right) .$ We have $\Delta =2\delta _{\max }-%
\mathcal{E}=se-e-d+1-(d+e)\left( e-l+m\right) .$

Putting $s=p-1=l(d+e)-1$ we get 
\begin{equation*}
\Delta =\left( d+e\right) (le-e+l-m)-d-2e+1.
\end{equation*}%
Since $m=\frac{s}{d}=l+e\frac{l}{d}-\frac{1}{d},$ $l,d\geq 2,$ the quantity 
\begin{equation*}
A:=le-e+l-m=e\left( l-l/d-1\right) +1/d
\end{equation*}%
is positive. If $A>1,$ then evidently $\Delta >0.$ But $A=1$ implies $%
le=(m-l)+e+1.$ Here $s=l(d+e)-1=md$ gives $le=(m-l)d+1,$ what yields $%
(m-l)(d-1)=e.$ Moreover, $m-l$ and $e$ are relatively prime, hence a
contradiction. $\qed$

\bigskip

\textit{Subcase} $(ii).$ Here $\nu _{0}=0.$ Because $\Delta
|_{|r|=s=r^{\prime }}=1-p^{\prime }$ may be strictly negative (the case is
not strict), we must check possibilities of hiding at $t=\infty $ and at the
singular point $t_{1}.$ Recall that $r^{\prime }>p^{\prime },2.$

\bigskip

3.15. \textbf{Lemma.} \emph{The possibility }$d=1$\emph{\ in (3.13) occurs
only in three cases (after normalizations):}%
\begin{eqnarray*}
\varphi &=&t^{2}(t-1),\;\;\;\psi =3t^{-1}-t^{-2}; \\
\varphi &=&t^{3}(t-1),\;\;\;\psi =2t^{-2}-t^{-3}; \\
\varphi &=&t^{3}(t-1),\;\;\;\psi =2t^{-2}+t^{-3}.
\end{eqnarray*}%
\emph{In the first two cases the (unique) double point becomes hidden at the
singularity }$t_{1};$\emph{\ in the third case the double point becomes
hidden at }$t=\infty .$\emph{\ These cases generate the series (g), (h) and
(i) of Main Theorem.\medskip }

\textit{Proof}. Assume that the double points are hidden at $t_{1}$; so $%
n_{1}=2\geq p^{\prime }.$ Then $2\delta _{t_{1}}=\mathcal{E}=2(s-\left\vert
q\right\vert )=2(p-\left\vert q\right\vert -1),$ $2\delta _{\max
}=s(p-\left\vert q\right\vert -1)-p^{\prime }+1$ (recall that $\left\vert
r\right\vert =r^{\prime }=s)$ and%
\begin{equation*}
\Delta =(s-2)(p-\left\vert q\right\vert -1)+1-p^{\prime }.
\end{equation*}

If $s=2,$ then $p=3,$ $\left\vert q\right\vert =1,$ $p^{\prime }=1$ and we
have the first of our cases.

Since $p-\left\vert q\right\vert -1\geq p^{\prime }-1,$ the other
possibility is $s=3.$ Then $\Delta =p-\left\vert q\right\vert -p^{\prime
}=0, $ where $p=4$ and $\left\vert q\right\vert \leq 2.$ Thus $p^{\prime
}=\left\vert q\right\vert =2$ and we get the second case from the lemma.

Let the double points be hidden at $t=\infty .$ So $p^{\prime }\geq 2.$ Then 
$2\delta _{\infty }\leq p^{\prime }\left( s-\left\vert q\right\vert \right)
=p^{\prime }(p-\left\vert q\right\vert -1)$ and 
\begin{equation*}
\Delta =(s-p^{\prime })(p-\left\vert q\right\vert -1)+1-p^{\prime }.
\end{equation*}%
It follows that $s=p^{\prime }+1$ and $\Delta =p-\left\vert q\right\vert
-p^{\prime }=0.$ Thus $p=s+1=p^{\prime }+2$ and $\left\vert q\right\vert =2,$
what implies $p^{\prime }=2,$ $p=4,$ $\left\vert r\right\vert =s=3.$ This is
the third of the cases from the lemma.

Finally we note that the polynomial curve $\left( \varphi ,\varphi \psi
\right) $ can be reduced either to a quasi-homogeneous curve (in the first
two cases) or a straight line (in the third case). $\qed$

\bigskip

Now Subcase B.2 is complete.

\bigskip

\textbf{B.3.} Recall that we have $p^{\prime }=r^{\prime }=1$ and the double
points are hidden at $t_{1}$ with $n_{1}=2.$ Here $\mathcal{E}=2\left(
s-\left\vert q\right\vert \right) $ and $\Delta =\left( \left\vert
r\right\vert -2\right) \left( s-\left\vert q\right\vert \right) +s-1>0$ (as $%
\left\vert r\right\vert \geq 2$ and $s>\left\vert q\right\vert \geq 1).$
There are no such annuli.

\bigskip

\textbf{3.III. The case }$p\geq \left\vert r\right\vert +2.$ Denote $%
K=p-\left\vert r\right\vert ,$ $L=s-\left\vert q\right\vert .$ Thus $2\leq
K\leq L.$ Here $\sigma =\dim Curv=K+L-1.$

We distinguish two subcases:%
\begin{equation*}
\text{A. }\max (p^{\prime },r^{\prime })\geq \max n_{j};\;\;\;\text{B. }%
n=\max n_{i}>p^{\prime },r^{\prime }.
\end{equation*}%
\medskip

\textbf{Case A.} Assume that $p^{\prime }$ is dominating, $p^{\prime }\geq
r^{\prime },n_{j}.$ Then $\mathcal{E}$ is maximal when all the double points
are hidden at $t=\infty $, i.e. $\mathcal{E}=p^{\prime }(K+L-1).$ We have $%
2\delta _{\max }=(K-1)(L-1)+\left[ p(\left\vert q\right\vert
+L)-(p-K)\left\vert q\right\vert \right] -p^{\prime }-r^{\prime }+1$ and 
\begin{equation}
\Delta =\left[ K\left\vert q\right\vert +Lp-(K+L)p^{\prime }\right]
+(KL-K-L+2)-r^{\prime }.  \tag{3.16}
\end{equation}%
The necessary condition for $\Delta \leq 0$ is%
\begin{equation*}
p=\left\vert q\right\vert =p^{\prime }
\end{equation*}%
(then the first term above vanishes).

Note now that $K+L=s-\left\vert r\right\vert $ is divisible by $r^{\prime }.$
Therefore $r^{\prime }\leq K+L$ and (3.16) yields%
\begin{equation*}
\Delta \geq (K-2)(L-2)-2.
\end{equation*}%
It is positive when $K,L\geq 4$ or when $K=3,$ $L\geq 5.$ We are left with
the cases%
\begin{equation}
K=2\leq L\text{ and }K=3,\text{ }L=3,4.  \tag{3.17}
\end{equation}

For $r^{\prime }$ dominating we get 
\begin{equation}
\Delta =\left[ Ks+L\left\vert q\right\vert -(K+L)r^{\prime }\right]
+(KL-K-L+2)-p^{\prime },  \tag{3.18}
\end{equation}%
the necessary condition 
\begin{equation*}
s=\left\vert r\right\vert =r^{\prime }
\end{equation*}%
and the same possibilities (3.17).

Let us stick to the case $p=\left\vert q\right\vert =p^{\prime }>r^{\prime
}. $ The quantity $s-\left\vert r\right\vert =K+L$ can be equal $r^{\prime }$
or $2r^{\prime }$ or can be greater. But if $K+L=2r^{\prime },$ then $\Delta
\geq (K-\frac{3}{2})(L-\frac{3}{2})-\frac{1}{4}$ and the only possibility is 
$K=L=2.$ We see also that $K+L\leq 2r^{\prime }.$

We distinguish three possibilities:%
\begin{eqnarray*}
\text{A.1. }p &=&\left\vert q\right\vert ,\text{ }r^{\prime }=2,\text{ }%
K=L=2; \\
\text{A.2. }p &=&\left\vert q\right\vert ,\text{ }K+L=r^{\prime }\text{;} \\
\text{A.3. }s &=&\left\vert r\right\vert ,\text{ }K+L=p^{\prime }.
\end{eqnarray*}%
(Due to the possible change $x\longleftrightarrow y,$ $t\rightarrow 1/t$ we
can assume $p^{\prime }>r^{\prime }$ when $K=L).$

\bigskip

\textbf{A.1.} Here $\sigma =3$ and $2\delta _{\max }=\mathcal{E}=3p.$ After
some rescalings we have 
\begin{equation*}
\varphi =t^{p}(1+at^{-1}+t^{-2}),\;\;\;\psi
=t^{-p}(1+b_{1}t^{-1}+b_{2}t^{-2}),
\end{equation*}%
and 
\begin{equation}
\varphi \psi =1+\gamma _{1}t^{-1}+\ldots +\gamma _{4}t^{-4}.  \tag{3.19}
\end{equation}%
We shall calculate the number $2\delta _{\infty }$ more carefully. Note that
the coefficients $\gamma _{j}$ in (3.19) determine the Puiseux quantities at 
$t=\infty ,$ e.g. $C_{1}^{(\infty )}=\gamma _{1}=a+b_{1}.$ \bigskip

3.20. \textbf{Lemma}. \emph{If }$\gamma _{1}=\gamma _{2}=\gamma _{3}=0,$%
\emph{\ i.e. }$\varphi \psi =1+\gamma _{4}t^{-4}$\emph{, then we obtain item
(s) from Main Theorem}.\medskip

\textit{Proof}. If $p$ were odd the equations $\varphi (t^{\prime })=\varphi
(t),$ $t^{\prime }=-t,$ would have a nontrivial solution (which would give a
solution to $\xi (t^{\prime })=\xi (t)).$ Let $p$ be even. If $a=0$ then
also $b_{1}=0$ and the curve is non-primitive (depends on $t^{2}).$ So let $%
a=-b_{1}\not=0$ (and $p$ even).

Each equation $\varphi (it)=\varphi (t)$ and $\psi (it)=\psi (t),$ i.e. for $%
t^{\prime }=it,$ has exactly one nontrivial solution $t_{1}$ and $t_{2}$
respectively. We must find the cases when these solutions are different,
i.e. $t_{1}\not=t_{2}.$ It can occur only when $\varphi (t_{1})=0$ or when $%
\psi (t_{2})=0.$ Indeed, only in this case the condition $\varphi \psi
(t^{\prime })=\varphi \psi (t)$ may not imply $\xi (t^{\prime })=\xi (t).$

But this means that the polynomials $1+at^{-1}+t^{-2}$ and $%
1+b_{1}t^{-1}+b_{2}t^{-2}$ have pairs of solutions $t_{1},$ $it_{1}$ and $%
t_{2},$ $it_{2}$ respectively. Now it is clear that these polynomials are $%
1\pm \sqrt{2}t^{-1}+t^{-2}.$ $\qed$

\bigskip

The next possibility is $\gamma _{1}=0\neq \gamma _{2}.$ Here we use the
estimate given in Lemma 2.32(c):%
\begin{equation}
2\delta _{\infty }\leq (2-1)(p-1)+(2-1)+2\nu _{\infty }^{\prime }=p+2\nu
_{\infty }^{\prime };  \tag{3.21}
\end{equation}%
here $2=p^{\prime \prime }=\gcd (p,2)$ and $\nu _{\infty }^{\prime }=\nu
_{\infty }-1$ is the `restricted' codimension of the degeneracy at $t=\infty
.$ However, we cannot apply (3.21) directly. The reason is that the
coefficient $p^{\prime \prime }=2$ before $\nu _{\infty }^{\prime }$ may be
smaller than $n_{1},$ which is $\leq 3.$ So, we should use rather the
estimate 
\begin{equation}
2\delta _{\infty }+2\delta _{0}+\sum 2\delta _{i}\leq p+3\cdot 2.  \tag{3.22}
\end{equation}%
But the latter is $<2\delta _{\max }=3p$ (for $p\geq 4).$

Therefore there remains the possibility $\gamma _{1}=\gamma _{2}=0\neq
\gamma _{3}.$ An analogous to (3.21) estimate gives $2\delta _{\infty }\leq
(3-1)(p-1)+(p^{\prime \prime }-1)+p^{\prime \prime }\nu _{\infty }^{\prime
}, $ where $p^{\prime \prime }=\gcd (p,3)=1$ or $=3.$ Similarly as in (3.22)
we get $2\delta _{\infty }+2\delta _{0}+\sum 2\delta _{i}\leq 2p+3<2\delta
_{\max }.$

It follows that there are no embedded annuli in the case A.1.

\bigskip

\textbf{A.2.} Recall that $p=\left\vert q\right\vert $ and $K+L=s-\left\vert
r\right\vert =r^{\prime }.$ Then $\sigma =K+L-1=r^{\prime }-1.$ We have 
\begin{equation}
\varphi =t^{p}(1+a_{1}t^{-1}+\ldots +a_{K}t^{-K}),\;\;\;\psi
=t^{-p}(1+b_{1}t^{-1}+\ldots +b_{_{L}}t^{-L})  \tag{3.23}
\end{equation}%
and 
\begin{equation}
\varphi \psi =1+\gamma _{1}t^{-1}+\ldots +\gamma _{r^{\prime }}t^{-r^{\prime
}}.  \tag{3.24}
\end{equation}

If $\gamma _{1}=\ldots =\gamma _{r^{\prime }-1}=0,$ then $2\delta _{\infty
}=(r^{\prime }-1)(p-1)+(p^{\prime \prime }-1),$ where $p^{\prime \prime
}=\gcd (p,r^{\prime })=\gcd (K,K+L)\leq K$ (as $p=\left\vert
r_{1}\right\vert r^{\prime }+K).$ Since $2\delta _{\max }=(r^{\prime
}-1)p+KL-2K-2L+2,$ we get 
\begin{equation}
\Delta \geq KL-2K-L+2=(K-1)(L-2).  \tag{3.25}
\end{equation}%
So we find $K=L=2.$ But this means that $\varphi \psi =1+\gamma _{4}t^{-4}$
and this possibility was considered in the point A.1.

Suppose that $\gamma _{1}=\ldots =\gamma _{l-1}=0\neq \gamma _{l},$ $%
l<r^{\prime }.$ Then $2\delta _{\infty }\leq (l-1)(p-1)+(p^{\prime \prime
}-1)+p^{\prime \prime }\nu _{\infty }^{\prime },$ where $p^{\prime \prime
}=\gcd (p,l)$ and $\nu _{\infty }^{\prime }=\nu _{\infty }-(l-1).$ Since $%
p^{\prime \prime }<r^{\prime },$ the other double points are hidden at $t=0$
and $2\delta _{0}+2\delta _{\infty }\leq (l-1)(p-1)+(p^{\prime \prime
}-1)+r^{\prime }(r^{\prime }-l),$ which is maximal for $l=r^{\prime }-1$ and
is bounded by 
\begin{equation*}
\mathcal{E}=(K+L-2)p+K+L-\varepsilon ,
\end{equation*}%
$\varepsilon =(r^{\prime }-1)-p^{\prime \prime }\geq 0.$ Thus, with $p\geq
2K+L,$ we have 
\begin{equation}
\Delta \geq p+KL-3K-3L+2+\varepsilon \geq (K-2)(L-1)+\varepsilon . 
\tag{3.26}
\end{equation}%
Therefore the case is strict and it should be $K=2,$ $p=2K+L=4+L$ and $%
\varepsilon =0.$ The latter implies that $p$ is a multiple of $r^{\prime
}-1=L+1.$ It follows that $L=2,$ $p=6,$ $r^{\prime }=\left\vert r\right\vert
=4$, $s=8$. Moreover, $\gamma _{1}=\gamma _{2}=0$ and in (3.24) and the
first Puiseux quantity $C_{1}^{(0)}$ at $t=0$ vanishes.

\bigskip

3.27. \textbf{Lemma.} \emph{After normalizations we have }$\varphi
=t^{4}(t^{2}+t+\frac{2}{3})$\emph{, }$\psi =t^{-8}(t^{2}-t+\frac{1}{3})$%
\emph{\ in (3.23). It is the exceptional case (t) of Main Theorem. \medskip }

\textit{Proof}. We know that $\gamma _{1}=a_{1}+b_{1}=0$ and $\gamma
_{2}=a_{2}+b_{2}-a_{1}^{2}=0,$ thus $b_{2}=a_{1}^{2}-a_{2}.$ Next, $\varphi
^{2}\psi =a_{2}^{2}b_{2}+C_{1}^{(0)}t+\ldots ,$ where $%
C_{1}^{(0)}=2a_{1}a_{2}b_{2}+a_{2}^{2}b_{1}=a_{1}a_{2}(2a_{1}^{2}-3a_{2})=0.$
Putting $a_{1}=1$ (normalization), we easily find the other constants. $\qed$

\bigskip

Therefore for $p^{\prime }$ dominating we have detected only two annuli.

\bigskip

\textbf{A.3.} If $r^{\prime }$ is dominating, then we obtain formulas
analogous as above. Equation (3.24) is replaced by $\varphi \psi =\gamma
_{0}+\gamma _{1}t+\ldots +\gamma _{p^{\prime }}t^{p^{\prime }},$ $p^{\prime
}=K+L.$

If $\gamma _{1}=\ldots =\gamma _{p^{\prime }-1}=0,$ then $2\delta
_{0}=(p^{\prime }-1)(s-1)+(s^{\prime \prime }-1),$ $s^{\prime \prime }=\gcd
(s,p^{\prime })=\gcd (L,K+L)\leq L$ (and $s^{\prime \prime }=L$ only when $%
K=L).$ Instead of (3.25) we get $\Delta \geq (K-2)(L-1)+\varepsilon ,$ where 
$\varepsilon =0$ for $K=L$ and $\varepsilon >0$ otherwise. So $K=L=2,$ but
this case is eliminated in the same fashion as for $p^{\prime }\geq
r^{\prime }.$ The inequality (3.26) is replaced with $\Delta \geq
(K-1)(L-2)+\epsilon ,$ $\epsilon =(p^{\prime }-1)-s^{\prime \prime }.$ Again
we get $K=L=2.$

\bigskip

\textbf{Case B.} When $n=n_{1}$ dominates over $p^{\prime },r^{\prime }$,
then $\mathcal{E}=n(\sigma +2-n),$ $\sigma =K+L-1.$ Since $K\leq L$ we put $%
n=K+1$ (i.e. maximal possible) and we get $\mathcal{E}=(K+1)L.$ Therefore 
\begin{eqnarray*}
\Delta &\geq &(p-2)L+(|q|-1)K-p^{\prime }-r^{\prime }+2 \\
&\geq &(|r|L-r^{\prime })+[(2(|q|-1)-p^{\prime }+2]\geq p^{\prime
}+r^{\prime }>0.
\end{eqnarray*}

Curves of Type $\binom{-}{-}$ have been investigated completely. $\qed$

\bigskip

\section{Annuli of Type $\binom{-}{+}$}

Recall that here $-p\leq r<0<q,s$ and $\frac{q}{p}\notin \mathbb{Z}.$ It
implies that $p_{1}\geq 2$ and 
\begin{equation}
\sigma =\dim Curv=p-\left\vert r\right\vert +q+s-2-k,\;\;\;k=[q/p]. 
\tag{4.1}
\end{equation}

If $p=\left\vert r\right\vert ,$ i.e. $\varphi =t^{p},$ then we apply Lemma
3.2. Therefore in the sequel we assume%
\begin{equation*}
p\geq \left\vert r\right\vert +1.
\end{equation*}%
\medskip

\textbf{4.I. When double points escape to infinity.} In this subsection we
assume that either $n_{\infty }=p^{\prime }$ or $n_{0}=r^{\prime }$ is
dominating over the orders $n_{i}$ of eventual singular points. We shall
show that this cannot occur.

The cases%
\begin{equation*}
\text{A. }p^{\prime }\geq r^{\prime },n_{i};\;\;\;\;\;\text{B. }r^{\prime
}\geq p^{\prime },n_{i}\text{ and }r^{\prime }\not=p^{\prime }
\end{equation*}%
are treated separately.

\bigskip

\textbf{Case A.} Since $p^{\prime }$ is dominating, we estimate $2\delta
_{\infty }+\sum 2\delta _{i}$ by 
\begin{equation}
\mathcal{E}=p^{\prime }\sigma =\left( p^{\prime }\right)
^{2}(p_{1}+q_{1})+p^{\prime }(s-\left\vert r\right\vert -2-k).  \tag{4.2}
\end{equation}%
We have also%
\begin{equation}
2\delta _{\max }=\left( p^{\prime }\right) ^{2}p_{1}q_{1}-(p+q)+(\left\vert
r\right\vert -s-\left\vert r\right\vert s)+2ps-(p^{\prime }+r^{\prime })+2. 
\tag{4.3}
\end{equation}%
We see then that the coefficient before $\left( p^{\prime }\right) ^{2}$ in $%
\Delta $ equals%
\begin{equation*}
(p_{1}-1)(q_{1}-1)-1.
\end{equation*}%
In general it is $\geq 0,$ but for $q_{1}=1$ it equals $-1$ (recall that $%
p_{1}>1).$ We distinguish then two subcases:%
\begin{equation*}
\text{A.1. }p/q\in \mathbb{Z};\;\;\;\text{A.2. }p_{1},q_{1}>1.
\end{equation*}%
\medskip

\textbf{A.1.} This is the first place, where we use the handsomeness
property (see Definition 2.44 and Proposition 2.45), which states that 
\begin{equation*}
s\geq p^{\prime }.
\end{equation*}%
Of course, $k=[\frac{q}{p}]=0.$ So we get 
\begin{equation*}
\Delta =-\left( p^{\prime }\right) ^{2}+2ps-p^{\prime }(p_{1}+s-\left\vert
r\right\vert )-\left\vert r\right\vert s+\left\vert r\right\vert
-s-r^{\prime }+2.
\end{equation*}%
This expression is increasing in $s:$ $\frac{\partial \Delta }{\partial s}%
=(p-\left\vert r\right\vert -1)+(p-p^{\prime })$ and $p>p^{\prime }.$ For $%
s=p^{\prime }$ we have $\Delta =2\left( p^{\prime }\right)
^{2}(p_{1}-1)-p^{\prime }(p_{1}+1)+(\left\vert r\right\vert -r^{\prime })+2.$
Using $p^{\prime }\geq 2$ we find $\Delta \geq p^{\prime
}(3p_{1}-5)+(\left\vert r\right\vert -r^{\prime })+2>0.$

\bigskip

\textbf{A.2.} Let (firstly) $1<q_{1}<p_{1};$ thus $q_{1}\geq 2,$ $p_{1}\geq
3 $ and $k=0.$ We rewrite the term $2ps$ in (4.3) in the form $p^{\prime
}(p_{1}s)+r^{\prime }(ps_{1}).$ Then we obtain%
\begin{equation}
\Delta =\left( p^{\prime }\right) ^{2}(p_{1}q_{1}-p_{1}-q_{1})+p^{\prime
}(p_{1}s-p_{1}-q_{1}-s+\left\vert r\right\vert +1)+r^{\prime
}(ps_{1}-\left\vert r\right\vert s_{1}-1)+2.  \tag{4.4}
\end{equation}%
It is increasing in $p_{1}$ and $q_{1}$ and for $p_{1}=3,$ $q_{1}=2$ equals $%
\Delta =p^{\prime }(p^{\prime }+2s+\left\vert r\right\vert -4)+r^{\prime }%
\left[ (p-\left\vert r\right\vert )s_{1}-1\right] +2>0.$

Let $1<p_{1}<q_{1}.$ Then $k\geq 1$ and $\Delta $ in (4.4) is increased by
at least $p^{\prime }.$ So analogous estimates give $\Delta >0.$

\bigskip

\textbf{Case B.} Now $r^{\prime }$ is dominating. Then with $\mathcal{E}%
=r^{\prime }(p-|r|+q+s-2)$ and $2\delta _{\max
}=(p-|r|-1)(q+s-1)+(ps_{1}+q|r_{1}|)r^{\prime }-p^{\prime }-r^{\prime }+1$
we have 
\begin{eqnarray*}
\Delta &=&(p-|r|-1)(q+s-1-r^{\prime })+[(p-|r|)s_{1}+(|r_{1}|-1)q]r^{\prime
}-p^{\prime }+1 \\
&\geq &r^{\prime }s_{1}-p^{\prime }+1\geq r^{\prime }-p^{\prime }+1>0.
\end{eqnarray*}%
\medskip

\textbf{4.5. Remark.} We can say that infinity is not a safe place for
hiding, at least for curves of Type $\binom{-}{+}.$ In particular, there are
no smooth handsome annuli of Type $\binom{-}{+}.$ The reader will see that
the same property (with the series (r) of Main Theorem as the only
exception) holds for handsome curves of Types $\binom{+}{+}$ and $\binom{-+}{%
+-}.$ On the other hand, annuli of Type $\binom{-}{-}$ can be characterized
as those, whose double points are hidden at infinity.

\bigskip

\textbf{4.II. When double points prefer finite singularities:} $\max
n_{i}>p^{\prime },r^{\prime }.$ Here the cases%
\begin{equation*}
\text{A. }p-\left\vert r\right\vert +1\leq q+s;\;\;\;\text{B. }q+s\leq p-|r|
\end{equation*}%
are treated separately.

\bigskip

\textbf{Case A.} Let $n=n_{1}$ be the maximal among the $x$-orders of the
singular points. We can assume $t_{1}=1.$ We have $n>p^{\prime },r^{\prime }$
and 
\begin{equation}
n\leq p-|r|+1  \tag{4.6}
\end{equation}%
(because $\varphi ^{\prime }(t)$ can have zero of order $\leq p-|r|).$ But
when there are more singular points, then $\sum_{i=1}^{N}(n_{i}-1)\leq p-|r|$
and $n_{1}\leq p-|r|$ (with equality for $N=2,$ $n_{2}=2).$

Nevertheless, with $2\delta _{j}\leq n_{j}\nu _{j},$ $2\delta _{\infty }\leq
n_{\infty }\nu _{\infty }$ and $\sum_{j=0}^{\infty }ext\nu _{j}=\sigma $
fixed the quantity $\mathcal{E}=\sum_{j=0}^{\infty }n_{j}\nu _{j}$ is
maximal when there is one singular point and 
\begin{equation*}
\mathcal{E}(n)=n(\sigma +2-n).
\end{equation*}%
$\mathcal{E}(n)$ achieves this this maximum for $n=p-|r|+1.$ Note however
that, when $p-|r|+1=q+s$ and $k=0,$ then $\mathcal{E}(p-|r|+1)=\mathcal{E}%
(p-|r|).$ Anyway the maximal value of $\mathcal{E}$ (used in below
calculations) is 
\begin{equation*}
\mathcal{E}=(p-|r|+1)(q+s-1-k),\;\;\;k=[q/p].
\end{equation*}

Therefore, using (2.11), we get%
\begin{eqnarray*}
\Delta &=&-2(q+s-1)+(ps+|r|q)-(p^{\prime }+r^{\prime })+1+k(p-|r|+1) \\
&=&(|r|-2)q+(p-2)(s-1)+(p-p^{\prime })+(1-r^{\prime })+k(p-|r|+1).
\end{eqnarray*}%
Here $p\geq |r|+1,$ $p\geq 2p^{\prime }$ and $k\geq 0.$

If $|r|\geq 2,$ then $\Delta \geq (s-1)+p^{\prime }+(1-r^{\prime
})=(s-r^{\prime })+p^{\prime }>0.$ Therefore 
\begin{equation}
|r|=r^{\prime }=1  \tag{4.7}
\end{equation}%
and $\varphi (t)$ can be treated as a general polynomial of degree $p.$

Now we write 
\begin{equation*}
p=p_{1}p^{\prime },\;\;\;q=(kp_{1}+q_{2})p^{\prime },\;\;\;0<q_{2}<p_{1}
\end{equation*}%
and we get 
\begin{equation*}
\Delta =(p_{1}-q_{2}-1)p^{\prime }+(p-2)(s-1).
\end{equation*}%
In order that $\Delta \leq 0$ we should have 
\begin{equation}
q_{2}=p_{1}-1,\text{ i.e. }q=(k+1)p-p^{\prime },  \tag{4.8}
\end{equation}%
and one of the two:%
\begin{equation*}
\text{A.1. }p=2;\;\;\;\;\;\;\;\text{A.2. }s=1.
\end{equation*}%
Moreover, the both cases are strict.

\bigskip

\textbf{A.1.} Here we can assume 
\begin{equation}
\varphi =t(t-1)  \tag{4.9}
\end{equation}%
and $q=2k+1$ (i.e. $p^{\prime }=q_{2}=1).$

But we know the form (4.9) of $\varphi ,$ we have met it in Subsection 3.II.
Here we can also apply the reverse tower transformation $\mathcal{T}:\psi
\rightarrow \left[ \psi (t)-\psi (1)\right] /\varphi (t)$ (see Definition
3.4). The resulting curve has smaller $q=-\mathrm{ord}_{t=\infty }\psi .$
After applying $\mathcal{T}$ several times we fall into Type $\binom{-}{-}.$
Then we use Lemma 3.10, which says that our annuli belong to the series (b)
of Main Theorem.

Since the case is strict, there is no room for any changes in the above
choices. Hence we have detected all possible annuli.

\bigskip

\textbf{A.2.} Recall that $|r|=s=1.$ Accordingly to the discussion at the
beginning of Case A we should consider separately the following
possibilities:

(i) $n=p;$

(ii) $n=p-1$ when $q+s=p-|r|+1,$ i.e. $n=q,$ $p=q+1.$

\bigskip

In the case (i) we can write 
\begin{equation}
\varphi =(t-1)^{n},\;\;\;\psi =(t-1)^{m}Q(t)/t.  \tag{4.10}
\end{equation}

If $m\geq n,$ then the change $(\varphi ,\psi )\rightarrow (\varphi ,\psi
/\varphi )$ gives also a $\mathbb{C}^{\ast }$-embedding, but with smaller $%
q. $ Maybe we fall into Case B: $p-|r|+1>q+s.$

Let $m<n.$ Since the case is strict, we should have $2\delta _{1}=\mu =n\nu
=\mu _{\min }+n^{\prime }\nu ^{\prime }$, where $\nu ^{\prime }=\deg Q-k$
(by Proposition 2.52). By Lemma 2.27 it should be $m=n^{\prime }$, $\nu
^{\prime }=0$ and $\mu =\mu _{\min }=n(m-1)\leq p(p-2)$ (as $m\leq p-1).$ It
is smaller than $2\delta _{\max }=(p-2)q+(p+q)-p^{\prime }\geq
p^{2}-p-p^{\prime }+1$ (as $q\geq p-1).$ So $n\not=p.$

\bigskip

Consider the case (ii)$.$ Here $\varphi =(t-1)^{q}(t+\alpha )$, $\psi
=(t-1)^{m}Q(t)/t$, $2\delta _{\max }=q^{2}+q$ and $\mu =2\delta _{1}\leq 
\mathcal{E}=n\nu =q(q+1).$ Since the case is strict, we have $\mu =n\nu =\mu
_{\min }+n^{\prime }\nu ^{\prime }$ and this holds in three cases (see Lemma
2.27): ($\alpha $) $m=n^{\prime }<q,$ $\nu ^{\prime }=0$ (but here $\mu =\mu
_{\min }=q(n^{\prime }-1)<2\delta _{\max )},$ or ($\beta $) $m=q+1$ (but $%
\mu =(q-1)q<2\delta _{\max }),$ or finally ($\gamma $) $m=q,$ i.e. 
\begin{equation*}
\varphi =(t-1)^{q}(t+\alpha ),\;\;\;\psi =(t-1)^{q}(t+\beta )/t.
\end{equation*}%
Note that for $q=2$ there can exist another singular point $t_{2}.$

\medskip

4.11. \textbf{Lemma.} \emph{We have }$q=2,$\emph{\ }$\alpha =1,$\emph{\ }$%
\beta =-1/3.$\emph{\ This case can be transformed to a subcase of the case
(p) for }$k=0$\emph{\ of Main Theorem. }\medskip

\textit{Proof}. The case $C_{1}^{(1)}\not=0$ can be reduced to the case ($%
\beta )$ by a change $\psi -$const$\cdot \varphi .$ So $C_{1}^{(1)}=0$ and
the same change gives $\psi =$const$\cdot (t-1)^{q+2}/t$. So $n^{\prime }=2,$
$\nu ^{\prime }=0$ and $\mu =\mu _{\min }=(q-1)(q+1)+(2-1)=q^{2}<2\delta
_{\max }=q^{2}+q$ (by the strictness).

Therefore there exists another singular point and $q=2.$ The condition $\psi
^{\prime }(t_{2})=0$ gives the values of $\alpha $ and $\beta .$ We have $%
\varphi =(t-1)^{2}(t+1)$, $\psi =(t-1)^{2}(t-\frac{1}{3})/t$ and $\psi -%
\frac{1}{3}\varphi =-\frac{1}{3}(t-1)^{4}/t.$ After the change $t\rightarrow
1/t$ we get the pair $((t-1)^{4}t^{-3},$ $(t-1)^{2}(t-3)t^{-2})$. $\qed$

\bigskip

\textbf{Case B.} Here $q+s\leq p-|r|.$ We denote by $n$ (respectively $%
n_{j}) $ the $y$-order of the singular point $t_{1}=1$ (respectively of $%
t_{j});$ respectively $\nu $ and $\nu _{j}$ denote the $x$-codimensions of $%
t_{1}=1$ and $t_{j}.$ Of course, $n>p^{\prime },r^{\prime }$ and $k=[\frac{q%
}{p}]=0.$

\bigskip

4.12. \textbf{Lemma.} \emph{It should be}%
\begin{equation*}
|r|=s=r^{\prime }=1
\end{equation*}%
\emph{and }$\Delta \geq 1-p^{\prime }.$ \medskip

\textit{Proof}. The quantity $\mathcal{E}=n(p-\left\vert r\right\vert
+q+s-n) $ is maximal for $n=q+s$ (no other singularities). Thus $\mathcal{E}%
=\left( q+s\right) (p-\left\vert r\right\vert )$ and 
\begin{eqnarray*}
\Delta &=&-(p-|r|)-(q+s)+(ps+|r|q)-(p^{\prime }+r^{\prime })+2 \\
&=&(p-|r|)(s-1)+(q+s)(|r|-1)-(p^{\prime }-1)-(r^{\prime }-1).
\end{eqnarray*}%
It is clear that $\Delta >0$ for $|r|>1.$ For $|r|=r^{\prime }=1$ and $s>1$
we get $\Delta \geq p-p^{\prime }\geq p^{\prime }>0.$ $\qed$

\bigskip

With $|r|=s=1$ we have 
\begin{equation}
2\delta _{\max }=pq+p-q-p^{\prime }.  \tag{4.13}
\end{equation}%
However the case is not strict, since $\Delta =1-p^{\prime }$ be negative.
Therefore we shall consider separately the following situations:%
\begin{equation*}
\begin{array}{cc}
\text{B.1. }n=q,\nu =p; & \text{B.2. }n=q,\nu =p-1,\nu _{2}=1; \\ 
\text{B.3. }n=q,\nu =p-1,\nu _{\infty }=1; & \text{B.4. }n=q+1,\nu =p-1.%
\end{array}%
\end{equation*}%
The situations B.1, B.2, B.3 will turn out strict, what allows to omit other
situations (like $n<q).$

\bigskip

\textbf{B.1.} Here $\mathcal{E}=pq$ and 
\begin{equation}
\Delta =p-q-p^{\prime }=0  \tag{4.14}
\end{equation}%
(by (4.13)). We see that the case is strict. The Milnor number $\mu =2\delta
_{1}$ must be equal $n\nu =pq$ as well as $\mu _{\min }+n^{\prime }\nu
^{\prime },$ where $n^{\prime }$ and $\nu ^{\prime }$ are defined for curves
of the form 
\begin{equation}
\varphi =(t-1)^{m}P(t),\;\;\;\psi =(t-1)^{n}(t+\beta )/t,\;\;\;n=q. 
\tag{4.15}
\end{equation}%
From Lemma 2.27 we know that this may occur in two situations:

~~~(i) $m$ is a multiple of $n,$ (ii) $m=dn+n^{\prime }$ and $\mu =\mu
_{\min }.$

In the case (i) we have either $m=n=q,$ or $m=p=2q$ (due to (4.14)).

If $m=q,$ then we take $\varphi -$const$\cdot \psi =(t-1)^{q+l}\widetilde{P}%
(t)/t$, $\deg \widetilde{P}=p^{\prime }+1-l.$ It gives $\mu \leq
(q-1)(q+l-1)+(p^{\prime \prime }-1)+p^{\prime \prime }(p^{\prime }+1-l),$
where $p^{\prime \prime }=\gcd (q,l)\leq l.$ The latter expression is
non-decreasing in $l$ for $l\leq p^{\prime }=p-q.$ For $l=p^{\prime }$ we
get $p^{\prime \prime }=p^{\prime }$ and $\mu \leq pq-(p+q-2p^{\prime
})<pq=2\delta _{\max }.$ For $l=p^{\prime }+1$ we have $\mu =(q-1)p<2\delta
_{\max }.$

If $m=2q=p,$ then from (4.15) we find that the first Puiseux constant $%
C_{1}^{(1)}=-\beta /(\beta +1)\neq 0$ and $\mu =(p-1)(q-1)+(q-1)<pq.$

In the case (ii) with $d=0,$ we have $n^{\prime }=m$ and $\mu =\mu _{\min
}=q(m-1)<pq.$ If $d\geq 1,$ i.e. $m=dq+n^{\prime }\leq p,$ then $\mu =\mu
_{\min }=q(dq+n^{\prime }-2)<pq.$

\bigskip

\textbf{B.2.} Here $\mathcal{E}=q(p-1)+2$ and 
\begin{equation*}
\Delta =p-p^{\prime }-2.
\end{equation*}%
Since $p\geq q+2$ (case B) and $n=q\geq 2,$ we get $p=4,$ $q=p^{\prime }=2$, 
$\nu =3$ and $\nu _{2}=1.$ The case is strict and we should have $\mu
=q(p-1).$ We can assume the form (4.15) with $m\geq 2.$

The case $m=p=4$ was discussed in the point B.1. The assumption $m=3$ yields 
$\mu =2<q(p-1).$

Therefore, with $\tau =t-1,$ we can assume 
\begin{equation}
\varphi =\tau ^{2}(1+\alpha \tau +\beta \tau ^{2}),\;\;\;\psi =\tau ^{2}%
\frac{1+\gamma \tau }{1+\tau }=\tau ^{2}+(\gamma -1)\tau ^{3}-(\gamma
-1)\tau ^{4}+\ldots  \tag{4.16}
\end{equation}%
This curve should have vanishing two essential Puiseux quantities at $\tau
=0 $ and an additional singularity. Computation of $C_{1}^{(1)}$ and $%
C_{3}^{(1)}$ gives 
\begin{equation}
\alpha =\gamma -1,\;\;\;\beta =1/2-\gamma .  \tag{4.17}
\end{equation}

Now $\varphi ^{\prime }=\tau M(\tau ),$ $\psi ^{\prime }=\tau N(\tau
)(1+\tau )^{-2},$ where $M=2+(3\gamma -3)\tau +(2-4\gamma )\tau ^{2}$ and $%
N=2+(3\gamma -1)\tau -2\gamma \tau ^{2}.$ The singular point $\tau _{2}$ is
determined by the nonzero zero of the polynomial $M-N,$ i.e. $\tau _{2}=%
\frac{2}{1-3\gamma }.$ But the condition $M(\tau _{2})=0$ implies $\gamma
=1. $ It gives a contradiction (see (4.16)).

\bigskip

\textbf{B.3.} Here $\mathcal{E}=q(p-1)+p^{\prime },$ $p^{\prime }\geq 2$ and 
\begin{equation*}
\Delta =p-2p^{\prime }.
\end{equation*}%
So $p=2p^{\prime }$ (the case is strict) and $n=q=p^{\prime },$ what
contradicts the assumption $n>p^{\prime }.$

\bigskip

\textbf{B.4.} Here 
\begin{equation*}
\varphi =(t-1)^{m}P(t),\;\;\;\psi =(t-1)^{n}/t,\;\;\;n=q+1.
\end{equation*}

We distinguish two possibilities:

\begin{equation*}
\text{(i) }q=1,\;\;\;\text{(ii) }q\geq 2.
\end{equation*}

(i) We make the change $\psi \rightarrow \psi +2=t+1/t.$

\bigskip

4.18. \textbf{Lemma.} \emph{If }$\psi =t+1/t,$\emph{\ then }$\varphi (t)$%
\emph{\ is a polynomial uniquely defined by the condition}%
\begin{equation}
\varphi (t)-\varphi (1/t)=(t-1)^{2l+1}(t+1)^{2u+1}t^{-l+u+1},  \tag{4.19}
\end{equation}%
\emph{\ i.e. we have the series (j) of Main Theorem}. \medskip

\textit{Proof}. The condition $\psi (t^{\prime })=\psi (t)$ (for a double
point) reads as $t^{\prime }=1/t.$ On the other hand, the only possible
singular points are $t=\pm 1.$ Therefore the only solutions to the equation $%
\varphi (t)=\varphi (1/t)$ are $t=1$ and $t=-1.$ From this the condition
(4.19) follows.

Let $t^{j}-t^{-j}=(t-1/t)P_{j-1}(\psi ),$ where $P_{i}$ are suitable
Chebyshev polynomials of degree $i.$ With $\varphi =\sum_{j=1}^{p}a_{j}t^{j}$
the condition (4.19) reads as $\sum a_{j}P_{j-1}(\psi )=(\psi -2)^{l}(\psi
+2)^{u}.$ Since the Chebyshev polynomials form a basis in the space of
polynomials of degree $\leq p-1$, we find that the latter equation for $%
a_{1},\ldots ,a_{p}$ has unique solution for fixed $l$ and $u.$ $\qed$

\bigskip

(ii) $n\geq 3.$ If $m\geq n$ in $\varphi =(t-1)^{m}P(t)$, $\psi =(t-1)^{n}/t$%
, the the change $\left( \varphi ,\psi \right) \rightarrow \left( \varphi
/\psi ,\psi \right) $ diminishes $p=\deg \varphi ;$ here we may fall into
the case A ($p<q+2),$ but we are not afraid of it.

Let $m<n.$ We use the estimate $2\delta _{1}\leq \mathcal{E}%
=mn-m-n+n^{\prime }(\deg P+1),$ $\deg P=p-m$ (for $n^{\prime }>1).$ The
quantity $\mathcal{E}=\mathcal{E}(m)$ becomes maximal for $m=n-n^{\prime },$
i.e. $\mathcal{E}=n^{2}-nn^{\prime }-2n+n^{\prime }(p-n+n^{\prime }+2).$
Now, for $\Delta $ we fix $n,n^{\prime }$ and vary $p;$ we should take it
minimal $p=n+1.$ Then $\mathcal{E}=n^{2}-nn^{\prime }+\left( n^{\prime
}\right) ^{2}-2n+3n^{\prime },$ $2\delta _{\max }=n^{2}+1-p^{\prime }$ and%
\begin{equation*}
\Delta =nn^{\prime }-\left( n^{\prime }\right) ^{2}+2n-3n^{\prime
}+1-p^{\prime }\geq \frac{n^{2}}{4}+\frac{n}{2}+1-p^{\prime }\geq \frac{n^{2}%
}{4}-\frac{n}{2}+2>0
\end{equation*}%
as $p^{\prime }\leq q=n-1.$

\bigskip

4.20. \textbf{Remark.} We could use the same arguments as in the point (ii)
for $n=2.$ Then we would arrive to the situation with $m=1,$ i.e. without
singularity at $t=1.$ But the point $t=-1$ may be singular.

\bigskip

Annuli of Type $\binom{-}{+}$ are complete. $\qed$

\section{Annuli of Type $\binom{+}{+}$}

These are most frequent annuli in the list in Main Theorem.

Recall that here%
\begin{equation*}
0<p<q,\;\;\;0<r<s,\;\;\;\min \left( q/p,s/r\right) \notin \mathbb{Z}%
,\;\;\;p^{\prime }\geq r^{\prime }
\end{equation*}%
(note that then $(p,r)\neq (1,1)$ and $p+r\geq 3)$ and 
\begin{equation}
\sigma =\dim Curv=p+r+q+s-3-k,\;\;\;k=\min \left( [q/p],[s/r]\right) . 
\tag{5.1}
\end{equation}

\bigskip

\textbf{5.I. Annuli with }$ps\neq rq.$ We begin with the case $p^{\prime
}\geq n_{i}$ ($x$-orders of singular points); recall that $p^{\prime }\geq
n_{0}=r^{\prime }$ (see 2.35).

\bigskip

5.2. \textbf{Lemma.} \emph{There are no annuli with }$p^{\prime }\geq n_{i}.$
\medskip

\textit{Proof}. The quantity $\mathcal{E}=p^{\prime }\nu _{\infty
}+\sum_{i=0}^{N}n_{i}\nu _{i}$ is maximal when $\nu _{\infty }=\sigma $ and $%
\nu _{i}=0$ for $i\geq 0.$ Thus $\mathcal{E}=\left( p^{\prime }\right)
^{2}(p_{1}+q_{1})+p^{\prime }(r+s-3-k)$, $2\delta _{\max }=\left( p^{\prime
}\right) ^{2}p_{1}q_{1}+p^{\prime }\left[ p_{1}s+rq_{1}+\left\vert
p_{1}s-rq_{1}\right\vert -p_{1}-q_{1}-1\right] +\left[ (r-1)(s-1)-(r^{\prime
}-1)\right] $ (where\linebreak $\left\vert p_{1}s+rq_{1}\right\vert \geq
r^{\prime })$ and hence%
\begin{equation}
\begin{array}{ccc}
\Delta & \geq & \left( p^{\prime }\right)
^{2}(p_{1}q_{1}-p_{1}-q_{1})+p^{\prime }\left[
(s-1)(p_{1}-1)+(r-1)(q_{1}-1)+k+r^{\prime }\right] \\ 
&  & +\left[ (r-1)(s-1)-(r^{\prime }-1)\right] .%
\end{array}
\tag{5.3}
\end{equation}

If $p_{1}>1,$ i.e. $\frac{q}{p}\notin \mathbb{Z},$ then the term $%
p_{1}q_{1}-p_{1}-q_{1}>0.$ The second term in (5.3) is $\geq p^{\prime
}(1+r^{\prime })$ and the third term is also positive (since $s>1).$
Therefore $\Delta >0.$

Let $p_{1}=1.$ The first term in (5.3) is $-\left( p^{\prime }\right) ^{2}.$
But, by the handsomeness property (see Definition 2.44 and Proposition 2.45)
we have $r\geq p^{\prime }.$ Then the second term in (5.3) satisfies%
\begin{equation*}
p^{\prime }\left[ (r-1)(q_{1}-1)+k+r^{\prime }\right] \geq p^{\prime }\left[
r-1+k+r^{\prime }\right] \geq \left( p^{\prime }\right) ^{2}+p^{\prime
}(k-1+r^{\prime }).
\end{equation*}%
We see that here also $\Delta >0.$ $\qed$

\bigskip

In the sequel we assume $n>p^{\prime }.$ The following simple inequality
will be often used.

\bigskip

5.4. \textbf{Lemma.} \emph{If }$ps\neq rq,$\emph{\ then}%
\begin{equation}
\mathcal{D}:=\left\vert ps-rq\right\vert -(p^{\prime }+r^{\prime }) +1\geq 0
\tag{5.5}
\end{equation}%
\emph{and equality holds only when }$p^{\prime }=1$\emph{\ or }$r^{\prime
}=1 $\emph{, and }$\left\vert p_{1}s_{1}-r_{1}q_{1}\right\vert =1.$ \medskip

Let us estimate the reserve $\Delta $ when it is minimal possible, i.e. with
one singular point and no degeneration at infinity. This minimal quantity
equals%
\begin{equation*}
\Delta =k(p+r)-2(q+s)+1+\mathcal{D}.
\end{equation*}%
We see that it may be negative, so the case is rather far from being strict.

The further analysis is divided into the following cases, depending on the
number $N$ of singular points:%
\begin{equation*}
\text{A. }N\geq 3;\;\;\;\text{B. }N=2,\;\;\;\text{C. }N=1.
\end{equation*}%
\medskip

\textbf{Case A.} Assume $t_{1}=1,t_{2},\ldots ,t_{N}$ are the singular
points with $n=n_{1}$ maximal. Then the quantity $\mathcal{E}=\sum n_{i}\nu
_{i}$ is maximal when $n_{2}=\ldots =n_{N}=2,$ $\nu _{2}=\ldots =\nu _{N}=1,$
$n=p+r+2-N$ and $\nu =(p+r+q+s-3-k)-(n-2)-(N-1)=q+s-2-k.$ So $\mathcal{E}%
=(p+r+2-N)(q+s-2-k)+2(N-1)$ and 
\begin{equation}
\Delta =(k+1)(p+r)+(N-3)(q+s)-(k+2)(N-2)+1-2(N-1)+\mathcal{D}.  \tag{5.6}
\end{equation}%
Here $\frac{\partial \Delta }{\partial N}=q+s-k-4\geq \left[ k(p+r)+2\right]
-k-4\geq k(p+r-1)-2\geq 2(k-1)\geq 0$ and the equality holds when $p+r=3,$ $%
q+s=5$ and $k=1.$ Next, 
\begin{equation*}
\Delta |_{N=3}=(k+1)(p+r-1)-4+\mathcal{D}.
\end{equation*}

We see that $\Delta \leq 0$ only when $p+r=3,$ $k=1$ and $\mathcal{D}=0;$
and the case is strict. The only possibility for this is (modulo the change $%
t\rightarrow 1/t)$%
\begin{equation*}
p=2,\;\;\;r=1,\;\;\;q=3,\;\;\;s=2,
\end{equation*}%
with $\sigma =4.$ After moving the corresponding critical values of $\varphi 
$ and $\psi $ (at $t=1)$ and extracting a const$\cdot \varphi $ from $\psi $
we can assume that 
\begin{equation}
\varphi =(t-1)^{2}(t+\alpha )t^{-1},\;\;\;\psi =(t-1)^{4}(t+\beta )t^{-2}. 
\tag{5.7}
\end{equation}%
We used the fact that the point $t=1$ has greatest $y$-codimension $\nu \geq
2$ (see the beginning of this point); so the $\mathrm{ord}_{t=1}\psi >3.$

\bigskip

5.8. \textbf{Lemma.} \emph{A curve of the form (5.7) represents a }$\mathbb{C%
}^{\ast }$\emph{-embedding only in the following situations:}

\begin{itemize}
\item $\alpha =2,$\emph{\ }$\beta =\frac{1}{2}$\emph{\ (item (u) of Main
Theorem);}

\item $\alpha =2(2\pm \sqrt{5}),$\emph{\ }$\beta =\frac{\alpha (2\alpha +1)}{%
3\alpha +2}$\emph{\ (item (v) of Main Theorem).}
\end{itemize}

\emph{In the first situation there is only one singularity }$\mathbf{A}_{8}$%
\emph{\ and in the second situation there are two singularities }$\mathbf{A}%
_{4}.$\emph{\ In particular, there cannot be three singularities.} \medskip

\textit{Proof}. We have $\varphi ^{\prime }=(t-1)A(t)t^{-2},$ $\psi ^{\prime
}=(t-1)^{3}B(t)t^{-3}$ with 
\begin{equation}
A=2t^{2}+\alpha t+\alpha ,\;\;\;B=3t^{2}+(2\beta +1)t+2\beta .  \tag{5.9}
\end{equation}

1. Suppose that there are three singular points $t_{1}=1,t_{2},t_{3}.$ Then $%
t_{2,3}$ should be common zeroes of the polynomials $A$ and $B.$ It follows $%
\frac{2\beta +1}{\alpha }=\frac{2\beta }{\alpha },$ what is impossible.

2. Suppose that there is singularity $\mathbf{A}_{6}$ at $t=1$ and a
singularity $\mathbf{A}_{2}$ at $t_{2}.$ With the notations $\tau =t-1$, $%
\alpha ^{\prime }=1/(\alpha +1),$ $\beta ^{\prime }=1/(\beta +1)$ we have 
\begin{equation}
\begin{array}{ccc}
\alpha ^{\prime }\varphi & = & \tau ^{2}\left[ 1+(\alpha ^{\prime }-1)\tau
-(\alpha ^{\prime }-1)\tau ^{2}+(\alpha ^{\prime }-1)\tau ^{3}+\ldots \right]
, \\ 
\beta ^{\prime }\psi & = & \tau ^{4}\left[ 1+(\beta ^{\prime }-2)\tau
-(2\beta ^{\prime }-3)\tau ^{2}+(3\beta ^{\prime }-4)\tau ^{3}+\ldots \right]
.%
\end{array}
\tag{5.10}
\end{equation}%
We see that the quantity $C_{5}^{(1)}$ is proportional to $\beta ^{\prime
}-2\alpha ^{\prime },$ so $\beta ^{\prime }=2\alpha ^{\prime }$ and $\beta
=2\alpha +1.$ Now $A=2t^{2}+(2\beta +1)t+(2\beta +1)$ and the difference $%
B-A=t^{2}-1;$ thus $t_{2}=-1.$ But $A(-1)=2\neq 0.$

3. Assume the singularity $\mathbf{A}_{8}$ at $t=1.$ So $C_{5}^{(1)}=0$,
i.e. $\beta ^{\prime }=2\alpha ^{\prime }$ and $C_{7}^{(1)}=0.$ From (5.10)
we get 
\begin{equation*}
\chi :=\beta ^{\prime }\psi -\left( \alpha ^{\prime }\varphi \right)
^{2}=-\left( \alpha ^{\prime }\right) ^{2}\tau ^{6}(1-2\tau +\ldots )
\end{equation*}%
and 
\begin{equation*}
\left( \alpha ^{\prime }\right) ^{-2}\chi +\left( \alpha ^{\prime }\varphi
\right) ^{3}=\tau ^{7}(3\alpha ^{\prime }-1)+\ldots .
\end{equation*}%
Therefore $\alpha ^{\prime }=\frac{1}{3}$ and $\beta ^{\prime }=\frac{2}{3},$
what gives $\alpha =2$ and $\beta =\frac{1}{2}.$

4. Assume two singularities $\mathbf{A}_{4}.$ One $\mathbf{A}_{4}$
singularity at $t=1$ exists. Another should be at a point $t_{2}$ which we
are going to find.

Consider 
\begin{equation*}
\upsilon :=\frac{\psi }{\varphi ^{2}}=\frac{t+\beta }{(t+\alpha )^{2}}%
=(\alpha -\beta )\left[ \frac{1}{t+\alpha }-\frac{1}{2\alpha -2\beta }\right]
^{2}+const.
\end{equation*}%
It follows that $t_{2}$ is defined from the equation $\frac{1}{t+\alpha }=%
\frac{1}{2\alpha -2\beta },$ i.e. $t_{2}=\alpha -2\beta .$

Consider the double point equations $\varphi (t^{\prime })=\varphi (t),$ $%
\upsilon (t^{\prime })=\upsilon (t).$ Denoting $u=t^{\prime }+t,$ $%
v=tt^{\prime }$ we get 
\begin{equation*}
u+(\alpha -2)-\alpha /v=0,\;\;\;v+\beta u+\alpha (2\beta -\alpha )=0.
\end{equation*}%
They imply the equation 
\begin{equation*}
-\beta u^{2}+(\alpha ^{2}3\alpha \beta +2\beta )u+\alpha (\alpha
^{2}-2\alpha \beta -2\alpha +4\beta -1)=0,
\end{equation*}%
whose left-hand side should equal $-\beta (u-2\alpha +4\beta )^{2}$ (as the
double points outside $t=1$ are hidden at $t_{2}).$ So we get the equations 
\begin{equation*}
\alpha ^{2}-7\alpha \beta +8\beta ^{2}+2\beta =0,\;\;\;3\alpha ^{2}-10\alpha
\beta +8\beta ^{2}+\alpha =0.
\end{equation*}%
From their difference we obtain $\beta =\alpha (2\alpha +1)/(3\alpha +2)$
and $t_{2}=\alpha ^{2}/(3\alpha +2).$

Now the condition $A(t_{2})=0$ (see (5.9)) gives 
\begin{equation*}
0=-\alpha ^{3}+7\alpha ^{2}+12\alpha +4=-(\alpha +1)(\alpha ^{2}-8\alpha -4).
\end{equation*}%
Since we reject the solution $\alpha =-1$ (as $\mathrm{ord}_{t=1}\varphi
=2), $ we find two values $\alpha =2(2\pm \sqrt{5}).$

In Main Theorem only one value $\alpha =2(2+\sqrt{5})$ is given. In fact,
the case with $\alpha =2(2-\sqrt{5})$ is equivalent to the previous one. The
equivalence is achieved by normalizing $t_{2}$ to $1$ and by renaming $%
t_{1}. $ We omit these calculations. $\qed$

\bigskip

\textbf{Case B:} (two singular points). Here we distinguish the subcases:%
\begin{equation*}
\text{B.1. }n\leq p+r-1;\;\;\;\;\;\text{B.2. }n=p+r.
\end{equation*}

\bigskip

\textbf{B.1.} For $n=p+r-1$ we have $\mathcal{E}=(p+r-1)(q+s-1-k)+2$ (see
(5.6)) and 
\begin{equation*}
\Delta =k(p+r-1)-2+\mathcal{D}.
\end{equation*}%
We see that $k=1,$ $p+r=3$ and $\mathcal{D}=0$ (the case is strict). As in
Case A we conclude that $p=2,$ $r=1,$ $q=3,$ $s=2,$ $\nu =3.$ We arrive to
the form (5.7) to which Lemma 5.8 applies.

Finally, since the case with $n=p+r-1$ is strict, the possibilities with
smaller $n$ do not occur.

\bigskip

\textbf{B.2.} We have $\mathcal{E}=(p+r)(q+s-2-k)+2$ and 
\begin{equation}
\Delta =(k+1)(p+r)-(q+s)-1+\mathcal{D},  \tag{5.11}
\end{equation}%
where we recall that $k=\min \left( [\frac{q}{p}],[\frac{s}{r}]\right) .$
Assume that $\frac{s}{r}<\frac{q}{p}$ and put%
\begin{equation*}
q=(k+1)p+q_{2},\;\;\;s=(k+1)r-s_{2},
\end{equation*}%
where $s_{2}$ should be positive. Then $\Delta \geq
s_{2}-q_{2}+(rq_{2}+ps_{2})-p^{\prime }-r^{\prime }$ (see (5.5).

If $q_{2}<0,$ then from (5.11) we get $\Delta >0.$ For $q_{2}\geq 0$ we can
write 
\begin{equation*}
\Delta \geq q_{2}(r-1)+(s_{2}-r^{\prime })+(ps_{2}-p^{\prime }).
\end{equation*}%
It can be $\leq 0$ only if all the above summands vanish; this case is
strict. Thus $s_{2}=r^{\prime }=p_{1}=1.$ Since $r>1$ (as $\min \left( \frac{%
q}{p},\frac{s}{r}\right) \notin \mathbb{Z})$, also $q_{2}=0.$ Therefore%
\begin{equation}
q=(k+1)p,\;\;\;s=(k+1)r-1,\;\;\;\mathcal{D}=0.  \tag{5.12}
\end{equation}%
(Note that $p=1,$ $r=q=2,$ $s=3$ satisfies (5.12); modulo the change $%
t\rightarrow 1/t$ it is a case studied above.) The range of cases in (5.12)
is too big to get some definite conclusions. To diminish this set we should
study in detail the Milnor number $\mu =2\delta _{1}.$

We have 
\begin{equation*}
\varphi =(t-1)^{n}t^{-r},\;\;\;\psi =(t-1)^{m}Q(t)t^{-s},\;\;\;n=p+r.
\end{equation*}%
Since the case is strict, it should be $\mu =n\nu =\mu _{\min }+n^{\prime
}\nu ^{\prime }.$ It occurs only in two cases:

~~~(i) $m=ln$ and (ii) $m=ln+n^{\prime },$ $\nu ^{\prime }=0.$

If $l\leq k$ in (i), then we apply the change $y\rightarrow y-$const$\cdot
x^{l}.$ But when $l\geq k+1$ we have $q+s\geq (k+1)(p+r),$ in contradiction
to (5.12).

In the case (ii) Lemma 2.27 gives 
\begin{equation}
\mu =\mu _{\min }\leq m(n-1)-n/2.  \tag{5.13}
\end{equation}%
We split this case into three subcases: ($\alpha )$ $m\leq q+s-2,$ ($\beta )$
$m=q+s-1$ and ($\gamma )$ $m=q+s$.

($\alpha )$ For $m\leq q+s-2$ the inequality (5.13) yields $\mathcal{E}\leq
(q+s-2)(p+r-1)-\frac{1}{2}(p+r)+2$ and%
\begin{equation*}
\Delta \geq 3(p+r-2)/2>0.
\end{equation*}

($\beta )$ For $m=q+s-1$ the same calculations give $\Delta \geq \frac{1}{2}%
(p+r-4).$ So $p+r\leq 4$ and the case is strict (as $\Delta \geq -\frac{1}{2}
$ means $\Delta \geq 0).$ We have then the following possibilities:

\begin{itemize}
\item $p=1,$ $r=2,$ $q=k+1,$ $s=2k+1,$ $m=3k+1;$

\item $p=1,$ $r=3,$ $q=k+1,$ $s=3k+2,$ $m=4k+2;$

\item $p=r=2,$ $q=2k+2,$ $s=2k+1,$ $m=4k+2.$
\end{itemize}

\bigskip

5.14. \textbf{Lemma.} \emph{Among the above possibilities only the first two
realize embedded annuli of the form:}

\begin{itemize}
\item $\varphi =(t-1)^{3}t^{-2},$\emph{\ }$\psi =(t-1)^{3k+1}(t-4)t^{-2k-1}$%
\emph{\ (item (k) of Main Theorem);}

\item $\varphi =(t-1)^{4}t^{-3},$\emph{\ }$\psi =(t-1)^{4k+2}(t-3)t^{-3k-2}$%
\emph{\ (item (p) of Main Theorem).}
\end{itemize}

\medskip

\textit{Proof}. We determine the parameter $\beta $ in $\psi
=(t-1)^{m}(t+\beta )t^{-s}$ from the condition $\psi ^{\prime }(t_{2})=0,$
where $t_{2}=-r/p$ is the second zero of $\varphi ^{\prime }(t).$ Here we
can assume $k=0.$

Therefore $\psi =(t-1)(t+\beta )/t=t+(\beta -1)-\beta /t,$ $\psi ^{\prime
}=1+\beta /t^{2}$ and we get $\beta =-4$ in the first case.

Next, $\psi =(t-1)^{2}(t+\beta )/t^{2}=t+(\beta -2)+(1-2\beta )/t+\beta
/t^{2},$ $\psi ^{\prime }=1+(2\beta -1)/t^{2}-2\beta /t^{3}$ and hence $%
\beta =-3$ in the second case.

Finally, $\psi =(t-1)^{2}(t+\beta )/t,$ $\psi ^{\prime }=2t+(\beta -2)-\beta
/t^{2}$ and $\psi ^{\prime }(-1)=-4\neq 0$ in the third case. $\qed$

\bigskip

($\gamma )$ For $m=q+s$ we have the following result, which is analogous to
Lemma 4.3 from [BZI].

\bigskip

5.15. \textbf{Lemma.} \emph{If }$\varphi =(t-1)^{n}t^{-r}$\emph{, }$\psi
=(t-1)^{m}t^{-s},$\emph{\ then we have one of the cases (l), (m), (n) of
Main Theorem.} \medskip

\textit{Proof}. (This proof is different from the proof of Lemma 4.3 in
[BZI].) The assumption $ps\neq rq$ implies that $\nu _{0}=\nu _{\infty }=0,$
as can be checked by direct calculations, and there is only one singular
point $t=1$ with $C_{1}^{(1)}\neq 0.$ (Note that for $ps=rq$ the curve is
multiply covered.) Therefore we have the equality $\mu =2\delta _{\max },$
i.e.%
\begin{equation}
(p+r-1)(q+s-1)+(n^{\prime }-1)=(p+r-1)(q+s-1)+\mathcal{D}.  \tag{5.16}
\end{equation}

If $n^{\prime }=1,$ then $\mathcal{D}=0.$ By Lemma 5.4 we can assume that $%
p^{\prime }=1$ and $p_{1}s_{1}-q_{1}s_{1}=1.$ It is item (l) of Main Theorem
with $r^{\prime }\rightarrow pl,$ $r_{1}\rightarrow n,$ $s_{1}\rightarrow l,$
$p\rightarrow m-pn,$ $q\rightarrow k-pl.$

Let $n^{\prime }>1$ and assume $p^{\prime }\geq r^{\prime }.$ Since $%
r^{\prime }n^{\prime }|(ns-mr)$ and $ps-rq=ns-mr,$ we have $\left\vert
ps-rq\right\vert \geq r^{\prime }n^{\prime }.$ Then (5.16) gives $n^{\prime
}-1\leq r^{\prime }n^{\prime }-p^{\prime }-r^{\prime }+1,$ i.e. 
\begin{equation*}
(r^{\prime }-1)(n^{\prime }-1)\geq (p^{\prime }-1).
\end{equation*}%
Therefore either ($1)$ $n^{\prime }=2$ and $p^{\prime }=r^{\prime },$ or ($%
2) $ $p^{\prime }=r^{\prime }=1$ and $\left\vert ps-rq\right\vert =n^{\prime
}.$

In the case ($1)$ it should be $p^{\prime }=r^{\prime }=2$ (as $\gcd
(p^{\prime },r^{\prime })|n^{\prime }).$ We get the series (n) of Main
Theorem.

In the case ($2)$ we get the series (m) of Main Theorem with $n^{\prime
}\rightarrow l.$ $\qed$

\bigskip

\textbf{Case C:} (one singular point). Let us divide it into two subcases:%
\begin{equation*}
\text{C.1. }n\leq p+r-1,\;\;\;\text{C.2. }n=p+r.
\end{equation*}
\medskip

\textbf{C.1.} For $n=p+r-1$ we have $\mathcal{E}=(p+r-1)(q+s-k)$ and 
\begin{equation}
\Delta =(k-1)(p+r-1)+\mathcal{D}.  \tag{5.17}
\end{equation}%
Therefore it should be 
\begin{equation}
k=1,\;\;\;\mathcal{D}=0  \tag{5.18}
\end{equation}%
and the case is strict (what allows to omit the cases with $n<p+r-1).$ It
follows that $\mu =n\nu =\mu _{\min }+n^{\prime }\nu ^{\prime }.$ By Lemma
2.27 it holds in two situations:

\begin{itemize}
\item $m=ln+n^{\prime },$ $\nu ^{\prime }=0;$

\item $m=ln.$
\end{itemize}

In the first case $\mathcal{E}=\mu _{\min }\leq m(n-1)-\frac{n}{2}\leq
(q+s)(p+r-2)-\frac{1}{2}(p+r-1)$ and $\Delta \geq (q+s)-\frac{1}{2}(p+r)+%
\frac{1}{2}>0.$

In the second case it should be $l=2.$ Indeed, the case $l=1$ is reduced via 
$y\rightarrow y-$const$\cdot x$. Recall that $k=1$ (by (5.18)). For $l\geq 3$
we get $q+s\geq 3(p+r-1)\geq 2(p+r).$ But, with $ps<rq,$ we have $s\leq
2r-1, $ $q\geq 2p+1$ and $\mathcal{D}\geq (p+r)-(p^{\prime }+r^{\prime
})+1>0.$

So $m=2(p+r-1)\leq q+s<2(p+r).$ We see that there are two possibilities:

~~~(i) $q+s=2(p+r)-2,$

~~~(ii) $q+s=2(p+r)-1.$

\bigskip

\textit{Subcase (i)}: $q+s=2(p+r)-2.$ We introduce the notations%
\begin{equation*}
q=2p-1+s_{2},\;\;\;s=2r-1-s_{2},\;\;\;s_{2}\geq 0.
\end{equation*}%
We have $\mathcal{D}=(s_{2}+1)p+(s_{2}-1)r-(p^{\prime }+r^{\prime })+1.$ If $%
s_{2}>1,$ then $\mathcal{D}>0$ (contradiction with (5.18)); therefore $%
s_{2}=0$ or $1.$ \medskip

\textit{Let} $s_{2}=1.$ Thus $q=2p=2p^{\prime }$, $s=2r-2$ and $\mathcal{D}%
=p^{\prime }-r^{\prime }+1=0.$ Moreover, $r^{\prime }=2$ and hence $%
p^{\prime }=1=p.$ Therefore we get a curve of the form 
\begin{equation*}
\varphi =(t-1)^{r}(t+\alpha )t^{-r},\;\;\;\psi =(t-1)^{2r}t^{2-2r}.
\end{equation*}%
\medskip

5.19. \textbf{Lemma.} \emph{The above curve always has a self-intersection.}
\medskip

\textit{Proof}. The condition $C_{1}^{(1)}=0$ leads to $\alpha =0.$ $\qed$

\bigskip

\textit{Let} $s_{2}=0.$ Thus $q=2p-1,$ $s=2r-1,$ so $p^{\prime }=r^{\prime
}=1$ and $\mathcal{D}=\left\vert p-r\right\vert -1=0.$ We can assume that $%
p=r+1,$ i.e. 
\begin{equation*}
\varphi =(t-1)^{2r}(t+\alpha )t^{-r},\;\;\;\psi =(t-1)^{4r}t^{1-2r}.
\end{equation*}%
\medskip

5.20. \textbf{Lemma.} \emph{The above curve is a }$\mathbb{C}^{\ast }$\emph{%
-embedding only for }$\alpha =1.$\emph{\ It is item (o) of Main Theorem.}
\medskip

\textit{Proof}. It follows from the condition $C_{1}^{(1)}=0.$ $\qed$

\bigskip

\textit{Subcase (ii)}: $q+s=2(p+r)-1.$ We put 
\begin{equation*}
q=2p+s_{2},\;\;\;s=2r-1-s_{2},\;\;\;s_{2}\geq 0.
\end{equation*}%
We have $\mathcal{D}=s_{2}(p+r)-(p^{\prime }+r^{\prime })+p+1=0.$ Therefore $%
s_{2}=0$ and we have the curve%
\begin{equation*}
\varphi =(t-1)^{p+r-1}(t+\alpha )t^{-r},\;\;\;\psi =(t-1)^{2(p+r-1)}(t+\beta
)t^{1-2r},
\end{equation*}%
where $r\geq 2$ (as $s>r).$

\bigskip

5.21. \textbf{Lemma.} \emph{It should be }$p=1,$\emph{\ }$r=2,$\emph{\ }$%
\alpha =\frac{1}{2},$\emph{\ }$\beta =2.$\emph{\ After the change }$%
t\rightarrow 1/t$\emph{\ it is the exceptional case (u) of Main Theorem.}
\medskip

\textit{Proof}. Since the case is strict, we should have the equality $\mu
_{\min }+n^{\prime }\nu ^{\prime }=n\nu ,$ where $n^{\prime }=n=p+r-1$ and $%
\nu ^{\prime }=2$ (there are two parameters). By Lemma 2.26 either $3|n$ and 
$C_{1}^{(1)}=C_{2}^{(1)}=0$ or $n=2$ and $C_{1}^{(1)}=C_{3}^{(1)}=0.$

Let $3|n.$ With $\tau =t-1,$ $\alpha ^{\prime }=1(\alpha +1),$ $\beta
^{\prime }=1/(\beta +1)$ we have 
\begin{eqnarray*}
\alpha ^{\prime }\varphi &=&\tau ^{n}\frac{1+\alpha ^{\prime }\tau }{(1+\tau
)^{r}}=\tau ^{n}\left[ 1+(\alpha ^{\prime }-r)\tau +\ldots \right] , \\
\beta ^{\prime }\psi &=&\tau ^{2n}\frac{1+\beta ^{\prime }\tau }{(1+\tau
)^{2r-1}}=\tau ^{2n}\left[ 1+(\beta ^{\prime }-2r+1)\tau +(2r-1)(r-\beta
^{\prime })\tau ^{2}+\ldots \right] .
\end{eqnarray*}%
The condition $C_{1}^{(1)}=0$ gives 
\begin{equation*}
\beta ^{\prime }=2\alpha ^{\prime }-1.
\end{equation*}%
Next, $\left( \alpha ^{\prime }\varphi \right) ^{2}=\tau ^{2n}\left[
1+(2\alpha ^{\prime }-2r)\tau +(2r^{2}-4\alpha ^{\prime }r+\left( \alpha
^{\prime }\right) ^{2})\tau ^{2}+\ldots \right] $ and $C_{2}^{(1)}=0$
implies $\alpha ^{\prime }=1,$ i.e. $\alpha =0$ (a contradiction).

So assume that $n=2,$ i.e. $p=1,$ $r=2.$ Then $(\beta ^{\prime }\psi
)/(\alpha ^{\prime }\varphi )^{2}-1=$\linebreak $-(\alpha ^{\prime
}-1)^{2}\tau ^{2}(1-2\alpha ^{\prime }\tau +\ldots )$ and $C_{3}^{(1)}=$const%
$\cdot (2-3\alpha ^{\prime }).$ Thus $\alpha ^{\prime }=\frac{2}{3},$ $\beta
^{\prime }=\frac{1}{3},$ i.e. $\alpha =\frac{1}{2},$ $\beta =2.$ $\qed$

\bigskip

\textbf{C.2:} $n=p+r.$ We have $\mathcal{E}=(p+r)(q+s-1-k)$ and 
\begin{equation}
\Delta =k(p+r)-(q+s)+1+\mathcal{D}.  \tag{5.22}
\end{equation}%
We can also assume that $m\leq q+s$ in $\psi =(t-1)^{m}Q(t)t^{-s},$
otherwise we use Lemma 5.13.

Before detecting the cases with $\Delta \leq 0$ we establish one useful fact.

\bigskip

5.23. \textbf{Lemma.} \emph{We have }$q+s\leq (k+1)(p+r).$ \medskip

\textit{Proof}. Suppose the reverse inequality and $\frac{s}{r}<\frac{q}{p}.$
We put 
\begin{equation*}
q=(k+1)p+q_{2},\;\;\;s=(k+1)r-s_{2},\;\;\;q_{2}>s_{2}\geq 0.
\end{equation*}%
Here $r>1$ as $\frac{s}{r}\notin \mathbb{Z}.$ We get 
\begin{equation*}
\Delta =(r-1)(q_{2}-1)+(p+1)(s_{2}-1)-(p^{\prime }-1)-(r^{\prime }-1).
\end{equation*}%
If $s_{2}\geq 2,$ then $q_{2}\geq 3$ and $\Delta >0;$ so $s_{2}=1.$ Thus $%
r^{\prime }=1$ and $\Delta =(r-1)(q_{2}-1)-(p^{\prime }-1)\leq 0$ iff $r=2$
and $q_{2}=p^{\prime },$ i.e. $q=(k+1)p+p^{\prime },$ $r=2,$ $s=2k+1,$ and
the case is strict.

By the strictness we should have $n\nu =\mu _{\min }+n^{\prime }\nu ^{\prime
}$ and either $m=ln+n^{\prime },$ $\nu ^{\prime }=0$ or $m=ln.$ But in the
first case $\mu \leq m(n-1)-\frac{n}{2}\leq (q+s-1)(p+r-1)-\frac{1}{2}(p+r)$
and $\Delta >0.$

In the second case it should be $m=(k+1)n$ (as $q+s=(k+1)(p+r)+(p^{\prime
}-1)<(k+2)(p+r)$ and the cases $l\leq k$ are destroyed via the changes $%
y\rightarrow y+$\textrm{const}$\cdot x^{l}).$ For $l=k+1$ we also apply the
change $y\rightarrow y+$\textrm{const}$\cdot x^{k+1}.$ After this the data $%
p,r,q,n$ remain unchanged, but $s\rightarrow (k+1)r=2(k+1).$ Altogether the
exponent $m$ changes to $\tilde{m}\leq q+s=m+(p^{\prime }-1)<(k+2)n.$ Hence $%
n$ does not divide $\tilde{m}$ and the previous arguments show that $\Delta
>0.$ $\qed$

\bigskip

Having the inequality from Lemma 5.23 we proceed further. Instead of the
bound $\mu \leq n\nu $ we use $\mu \leq \mu _{\min }+n^{\prime }\nu ^{\prime
},$ where 
\begin{equation*}
\nu ^{\prime }\leq \deg Q-k+[(m-1)/n]\leq q+s-m
\end{equation*}%
for $\psi =(t-1)^{m}Q(t)t^{-s}$ (by Proposition 2.52). Instead of (5.22) we
get 
\begin{equation}
\Delta \geq (n-n^{\prime }-1)(q+s-m)-(n^{\prime }-1)+\mathcal{D}.  \tag{5.24}
\end{equation}

Supposition $n|m$ leads to $(k+1)(p+r)\leq m\leq q+s$ (since the cases $%
\frac{m}{n}\leq k$ can be destroyed) but we assumed that $m<q+s.$ Thus $%
n\geq 2n^{\prime }$ and (5.24) gives 
\begin{equation*}
\Delta \geq (n^{\prime }-1)(q+s-m-1)+(n-2n^{\prime })+\mathcal{D}.
\end{equation*}%
Therefore the case is strict, $\mathcal{D}=0,$ $n=2n^{\prime }$ and either $%
n^{\prime }=1,$ or $m=q+s-1.$ But $n^{\prime }=1$ implies $p=r=1$ and $\frac{%
q}{p},\frac{s}{r}\in \mathbb{Z}$, in contradiction with (2.35). Let 
\begin{equation*}
m=q+s-1.
\end{equation*}%
The condition $\mathcal{D}=0$ implies that either $p^{\prime }=1$ or $%
r^{\prime }=1;$ assume $r^{\prime }=1.$ Then $\left\vert ps-rq\right\vert
=p^{\prime }.$

Since $n^{\prime }=\frac{1}{2}(p+r)$ and $m=ln^{\prime },$ we have $q+s=l%
\frac{p+r}{2}+1.$ Because $k(p+r)+1<q+s\leq (k+1)(p+r),$ we get $l=2k+1.$ So 
\begin{equation*}
p+r=2n^{\prime },\;\;\;q+s=(2k+1)n^{\prime }+1.
\end{equation*}%
Since $ps-rq=p(q+s)-q(p+r),$ we find 
\begin{equation}
ps-rq=n^{\prime }\left[ (2k+1)p-2q\right] +p;  \tag{5.25}
\end{equation}%
it equals $\pm p^{\prime }.$

We distinguish three subcases:%
\begin{equation*}
\text{(}i\text{) }q>(k+1)p;\;\;\;(\text{ii})\text{ }q=(k+1)p;\;\;\;(\text{iii%
})\text{ }q<(k+1)p.
\end{equation*}%
In the case (i$)$ we get $\left\vert ps-rq\right\vert >n^{\prime }p-p\geq
p^{\prime }$ (as $n^{\prime }\geq 2).$

In the case (ii$)$ we have $p=p^{\prime }$ and $p=\left\vert
ps-rq\right\vert =(n^{\prime }-1)p.$ Thus $n^{\prime }=2$, i.e. $p+r=4,$ and
there are two types of curves:%
\begin{equation}
\varphi =(t-1)^{4}t^{-3},\;\;\;\psi =(t-1)^{4k+2}(t+\beta )t^{-3k-2}, 
\tag{5.26}
\end{equation}%
\begin{equation}
\varphi =(t-1)^{4}t^{-2},\;\;\;\psi =(t-1)^{4k+2}(t+\beta )t^{-2k-1}. 
\tag{5.27}
\end{equation}%
\medskip

5.28. \textbf{Lemma.} \emph{In (5.26) we have }$\beta =1$\emph{\ (after the
change }$t\rightarrow 1/t$\emph{\ it is a subcase of the series (o) from
Main Theorem) and the curve (5.27) always has a self-intersection.} \medskip

\textit{Proof}. Apply $C_{1}^{(1)}=0.$ $\qed$

\bigskip

Assume the case (iii$),$ i.e. $q<(k+1)p.$ It is impossible that $(2k+1)p=2q,$
because then (5.25) implies $\left\vert ps-rq\right\vert =p>p^{\prime }.$ By
the same reason $2q$ cannot be smaller than $(2k+1)p.$ Therefore $%
2q-(2k+1)p\geq p^{\prime }$ and hence 
\begin{equation}
\left\vert ps-rq\right\vert \geq |n^{\prime }p^{\prime }-p|=\left\vert
p\left( \frac{p^{\prime }}{2}-1\right) +\frac{r}{2}p^{\prime }\right\vert . 
\tag{5.29}
\end{equation}%
Since $\left\vert ps-rq\right\vert =p^{\prime },$ then either 
\begin{equation*}
\text{(}\alpha \text{) }p^{\prime }=1,\;\;\;\;\text{ or (}\beta \text{) }%
p^{\prime }=2.
\end{equation*}

($\alpha $) If $p^{\prime }=1,$ then $1\geq \frac{\left\vert r-p\right\vert 
}{2}.$ We can then assume $p=r-2$ (the supposition $p=r$ leads to $p^{\prime
}=p=r=1).$ Moreover, it should be $2q-(2k+1)p=p^{\prime }=1.$ Thus $p=2d-1$
is odd, $r=2d+1$, $q=\frac{1}{2}\left[ (2k+1)p-1\right] =k(2d-1)+d$ and $%
s=k(2d+1)+k+1.$ The eventual curve is of the form%
\begin{equation*}
\varphi =(t-1)^{4d}t^{-2d-1},\;\;\;\psi =(t-1)^{2d(2k+1)}(t+\beta
)t^{-k(2d+1)-k-1}.
\end{equation*}%
\medskip

5.30. \textbf{Lemma.} \emph{We have }$\beta =1$\emph{, i.e. item (o) of Main
Theorem (after the change }$t\rightarrow 1/t).$ \medskip

\textit{Proof}. Calculate $C_{1}^{(1)}.$ $\qed$

\bigskip

($\beta $) Let $p^{\prime }=2.$ The above formulas give $\left\vert
ps-rq\right\vert =r=2$ (by (5.29)) and also $2q-(2k+1)p=2.$ It gives $%
p=2(2d-1),$ $q=(2k+1)(2d-1)+1=2\left[ k(2d-1)+d\right] ,$ $s=2k+1$ and $%
s=2k+1,$ i.e. 
\begin{equation*}
\varphi =(t-1)^{4d}t^{-2},\;\;\;\psi =(t-1)^{2d(2k+1)}(t+\beta )t^{-2k-1}.
\end{equation*}%
\medskip

5.31. \textbf{Lemma.} \emph{The above curve always has a self-intersection.}
\medskip

\textit{Proof}. We have $\varphi =\tau ^{4d}(1-2\tau +\ldots ),$ $\beta
^{\prime }\psi =\varphi ^{k}\cdot \tau ^{2d}\left[ 1+(\beta ^{\prime
}-1)\tau +\ldots \right] ,$ $\tau =t-1,$ $\beta ^{\prime }=1/(\beta +1).$ So 
$C_{1}^{(1)}=0$ means $\beta ^{\prime }=0.$ $\qed$

\bigskip

5.32. \textbf{Remark.} Above, i.e. in the case with maximal $n$ and with one
singular point, we sometimes used the quantity $\mathcal{E}=\mu _{\min
}+n^{\prime }\nu ^{\prime }.$ But when $p^{\prime }>n^{\prime }$ or $%
r^{\prime }>n^{\prime },$ the double points (not $\frac{1}{2}\mu _{\min }$
double points hidden at $t=0)$ may prefer to hide themselves at infinity.
Therefore we should consider separately the quantity 
\begin{equation*}
\mathcal{E}=\mu _{\min }+n^{\prime }\nu ^{\prime }+n_{0}\nu _{0}+n_{\infty
}\nu _{\infty },
\end{equation*}%
where $\nu ^{\prime }+\nu _{0}+\nu _{\infty }\leq \deg Q.$

Suppose that $p^{\prime }>n$ and $p^{\prime }\geq r^{\prime }.$ Then we get $%
\mathcal{E}\leq \mu _{\min }+p^{\prime }(q+s-m)$ and 
\begin{equation*}
\Delta \geq (p+r-1-p^{\prime })(q+s-m)-(n^{\prime }-1)+\mathcal{D}.
\end{equation*}%
If $p=p^{\prime },$ then we use the handsomeness property $r\geq p;$ so $%
\Delta \geq (p^{\prime }-1)-(n^{\prime }-1)+\mathcal{D}>0.$ If $p\geq
2p^{\prime },$ then also $\Delta >0.$ No new cases are found.

\bigskip

Therefore the type $\binom{+}{+}$ with $ps\neq rq$ is complete. $\qed$

\bigskip

\textbf{5.II. The case }$ps=rq.$ By Remark 2.33 we have 
\begin{equation*}
p=p_{1}p^{\prime },\;\;\;r=p_{1}r^{\prime },\;\;\;q=q_{1}p^{\prime
},\;\;\;s=q_{1}r^{\prime },\;\;\;1<q/p\notin \mathbb{Z}.
\end{equation*}%
Thus $p_{1}\geq 2,$ $q_{1}\geq p_{1}+1\geq 3.$

We introduce also the notation%
\begin{equation*}
d^{\prime }=p^{\prime }+r^{\prime }.
\end{equation*}%
Therefore 
\begin{equation}
2\delta _{\max }=(p_{1}d^{\prime }-1)(q_{1}d^{\prime }-1)-(d^{\prime }-1). 
\tag{5.33}
\end{equation}

By Proposition 2.28(b) twice the number of hidden double points is estimated
by 
\begin{equation}
\mathcal{E}=\sum_{i=1}^{N}n_{i}\nu _{i}+d^{\prime }(\nu _{\inf }+1), 
\tag{5.34}
\end{equation}%
where 
\begin{equation}
\sum_{i=1}^{N}(n_{i}-2+\nu _{i})+\nu _{\inf }\leq \sigma
=(p_{1}+q_{1})d^{\prime }-3-k,\;\;\;k=\left[ q_{1}/p_{1}\right] .  \tag{5.35}
\end{equation}%
Moreover, when $\nu _{\inf }=0,1$ or when $d^{\prime }=2,$ then $d^{\prime
}(\nu _{\inf }+1)$ can be replaced by $d^{\prime }\nu _{\inf }.$

\bigskip

5.36. \textbf{Lemma.} \emph{In the case }%
\begin{equation*}
d^{\prime }\geq \max n_{i}.
\end{equation*}%
\emph{there are no such curves}.\medskip

\textit{Proof}. Indeed, by (5.34) and (5.35) we have $\mathcal{E}\leq
d^{\prime }(\sigma +1)=\left( d^{\prime }\right) ^{2}(p_{1}+q_{1})-d^{\prime
}(2+k)$ and $\Delta \geq \left( d^{\prime }\right)
^{2}(p_{1}q_{1}-p_{1}-q_{1})-d^{\prime }(p_{1}+q_{1}-1-k)+2.$ It is
increasing in $p_{1}$ and $q_{1}.$ For $p_{1}=2$ and $q_{1}=kp_{1}+1=2k+1$
we get 
\begin{equation*}
\Delta \geq \left( d^{\prime }\right) ^{2}(2k-1)-d^{\prime }k+2\geq
d^{\prime }(3k-2)+2>0.
\end{equation*}%
$\qed$\medskip

Now consider the case when 
\begin{equation*}
n=n_{1}>d^{\prime }
\end{equation*}%
and $n_{1}$ is maximal among $n_{i}.$ We divide the problem into two cases:%
\begin{equation*}
\text{A. }\nu _{\inf }\geq 2\text{ and }d^{\prime }\geq 3,\;\;\;\;\;\;\text{%
B. }\nu _{\inf }\leq 1\text{ or }d^{\prime }=2
\end{equation*}%
(this division is motivated by Lemma 2.32).

\medskip

\textbf{Case A.} Here $\mathcal{E}$ becomes maximal when $\nu _{\inf }=2$
and there is one singular point with maximal $n=p+r=p_{1}d^{\prime }$ (i.e. $%
\varphi =(t-1)^{n}t^{-r})$ and with $\nu =\sigma -(n-2)-\nu _{\inf
}=q_{1}d^{\prime }-3-k.$ Thus $\mathcal{E}=p_{1}d^{\prime }(q_{1}d^{\prime
}-3-k)+3d^{\prime }$ and 
\begin{equation*}
\Delta =(k+2)p_{1}d^{\prime }-q_{1}d^{\prime }-4d^{\prime }+2\geq
(p_{1}-3)d^{\prime }+2
\end{equation*}%
(since $q_{1}\leq (k+1)p_{1}-1).$ Of course, it should be $p_{1}=2,$ $%
q_{1}=2k+1$ and%
\begin{equation*}
\Delta =2-d^{\prime };
\end{equation*}%
we see that the case is not strict.

Assume $\psi =(t-1)^{m}Q(t)t^{-s}.$ If $m\geq n,$ then we make the change $%
\psi \rightarrow \psi /\varphi ;$ so assume that $m<n=2d^{\prime }.$ Here $%
\deg Q=(2k+1)d^{\prime }-m>(2k-1)d^{\prime }$ and $n^{\prime }\leq d^{\prime
}.$ We use $\mathcal{E}=\mu _{\min }+n^{\prime }\nu ^{\prime }+d^{\prime
}(\nu _{\inf }+1),$ which is maximal when $\nu ^{\prime }=0$ and $\nu _{\inf
}=\deg Q.$ Hence $\mathcal{E}=m(n-1)-(n-n^{\prime })+d^{\prime }\left[
(2k+1)d^{\prime }-m+1\right] ;$ it is increasing in $m.$ For $m=n-n^{\prime
} $ we get $\mathcal{E}=(n-n^{\prime })(n-d^{\prime }-2)+d^{\prime }\left[
(2k+1)d^{\prime }+1\right] .$ Since $n=2d^{\prime }$ and $n-n^{\prime
}<2d^{\prime },$ we have $\mathcal{E}<(2k+3)(d^{\prime })^{2}-3d^{\prime }$
and 
\begin{equation*}
\Delta >(2k-1)(d^{\prime })^{2}-(2k+1)d^{\prime }+2\underset{k=1}{\geq }%
(d^{\prime })^{2}-3d^{\prime }+2>0
\end{equation*}%
as $d^{\prime }\geq 3.$ Also when $n<p+r$ the same analysis leads to $\Delta
>0.$

\bigskip

\textbf{Case B:} $\nu _{\inf }\leq 1$ or $d^{\prime }=2.$ Recall that by
Lemma 2.32 we have $2\delta _{\inf }\leq d^{\prime }\nu _{\inf }.$ Since $%
d^{\prime }<n,$ then the double points would rather prefer finite
singularity $t_{1}=1.$ For fixed number $N$ of finite singularities $%
\mathcal{E}$ is maximal when $\nu _{\inf }=0$ and $2\delta _{t_{j}}=n_{j}\nu
_{j}=2\cdot 1$ for $j>1.$ In fact, the case with $d^{\prime }=2$ (and $\nu
_{\inf }=1)$ can be treated in the same way as an additional cusp point; in
that case we agree to regard $N$ as the number of singular points plus $1.$
Here 
\begin{equation*}
\mathcal{E}=n\left[ (p_{1}+q_{1})d^{\prime }-3-k-(n-2)-(N-1)\right] +2(N-1).
\end{equation*}

Two subcases are distinguished:%
\begin{equation*}
\text{B.1. }N\geq 2;\;\;\;\;\;\;\text{B.2. }N=1.
\end{equation*}

\bigskip

\textbf{B.1.} From the below estimates it will follow that $\Delta >0$ for $%
N>2;$ so we put $N=2.$ Then $\mathcal{E}$ is maximal for $n=p_{1}d^{\prime }$
(maximal), $\mathcal{E}=p_{1}d^{\prime }(q_{1}d^{\prime }-2-k)+2$ and%
\begin{equation*}
\Delta =\left[ (k+1)p_{1}-q_{1}-1\right] d^{\prime }.
\end{equation*}%
We see that $\Delta \leq 0$ iff $q_{1}=(k+1)p_{1}-1$ and the case is strict.

With more precise estimate of the Milnor number we should have $\mu =n\nu
=\mu _{\min }+n^{\prime }\nu ^{\prime }$. So, either $n|m$ (here we use $%
\psi \rightarrow \psi /\varphi )$ or $m=n^{\prime }<n$ (Lemma 2.27) and $%
m\leq q_{1}d^{\prime }-1.$ In the latter case (without assumption $%
p_{1}<q_{1})$ we have $\mu \leq m(n-1)-\frac{n}{2}\leq (q_{1}d^{\prime
}-1)(p_{1}d^{\prime }-1)-\frac{1}{2}p_{1}d^{\prime }$ and%
\begin{equation*}
\Delta \geq (p_{1}/2-1)d^{\prime }+1>0.
\end{equation*}

\bigskip

\textbf{B.2} (one singular point). We make additional division:%
\begin{equation*}
\text{(i) }n\leq p_{1}d^{\prime }-1,\;\;\;\;\;\;\text{(ii) }n=p_{1}d^{\prime
}.
\end{equation*}

\bigskip

In the \textit{subcase} (i$)$ (with $\nu _{\inf }=0$ and $%
n_{1}=p_{1}d^{\prime }-1)$ we have $\mathcal{E}=(p_{1}d^{\prime
}-1)(q_{1}d^{\prime }-k)$ and 
\begin{equation*}
\Delta =(k-1)(p_{1}d^{\prime }-1)-(d^{\prime }-1).
\end{equation*}%
We see that $\Delta \leq 0$ iff $k=1;$ here $\Delta =1-d^{\prime }<0$ and
the case is definitely not strict. Moreover, we can assume that $%
n=p_{1}d^{\prime }-1,$ as $\Delta >0$ otherwise.

We use more precise estimate $\mathcal{E}=(p_{1}d^{\prime
}-2)(q_{1}d^{\prime }-1)+(n^{\prime }-1)+n^{\prime },$ i.e. for maximal $%
\mathrm{ord}\psi _{t=1}=m=q_{1}d^{\prime }$. Here 
\begin{equation*}
\Delta =(q_{1}-1)d^{\prime }-2n^{\prime }+1.
\end{equation*}%
If $n^{\prime }\leq \frac{n}{2}=\frac{1}{2}(p_{1}d^{\prime }-1),$ then $%
\Delta \geq (q_{1}-p_{1}-1)d^{\prime }+\frac{3}{2}>0.$ So $n^{\prime }=n,$ $%
m=2n=2p_{1}d^{\prime }-2,$ $p^{\prime }=r^{\prime }=1$ and $\Delta
=1-d^{\prime }=-1.$ We have $p=r,$ $q=s=2p-1,$ i.e. 
\begin{equation*}
\varphi =(t-1)^{2p-1}(t+\alpha )t^{-p},\;\;\;\psi =(t-1)^{2(2p-1)}t^{1-2p}.
\end{equation*}%
\medskip

5.37. \textbf{Lemma.} \emph{The above curve is a }$\mathbb{C}^{\ast }$\emph{%
-embedding iff either }$\alpha =1$\emph{\ (item (q) of Main Theorem) or }$%
p=2,$\emph{\ }$\alpha =e^{\pm i\pi /3}$\emph{\ (item (r) of Main Theorem).}
\medskip

\textit{Proof}. The first statement follows directly from calculation of $%
C_{1}^{(1)}.$

The second case arises from the possibility of a double point hidden at
infinity, $2\delta _{\inf }=d^{\prime }\nu _{\inf }=2\cdot 1.$ We have $%
2\delta _{1}=\mu _{\min }=(2p-2)(4p-3)+(2p-2)=8p^{2}-12p+4$ and $2\delta
_{\max }=(2p-1)(4p-3)-1=8p^{2.}-10p+2.$ So it should be $p=2.$

But here we have the condition of equality of the leading terms of the
Puiseux expansions at $t=0$ and at $t=\infty .$ As $t\rightarrow \infty $ we
have $x\sim t^{p}=t^{2},$ $y\sim t^{3}\sim x^{3/2}.$ As $t\rightarrow 0$ we
have $x\sim (-\alpha )t^{-2},$ $y\sim t^{-3}$ and $y^{2}/x^{3}\sim (-\alpha
)^{-3}.$ Therefore $\alpha ^{3}=-1$ and $\alpha $ takes two values given in
the thesis of the lemma.

But the curve with $\alpha =e^{-i\pi /3}$ is equivalent to the curve with $%
\alpha =e^{i\pi /3}$ via the change $t\rightarrow 1/t$ and a normalization
of $x.$ $\qed$

\bigskip

5.38. \textbf{Remark}. Items (q) and (r) of Main Theorem present series
which include the curves from Lemma 5.37 as particular cases. These series
are produced using the tower transformations $(\varphi ,\psi )\rightarrow
(\varphi \psi ^{j},\psi ).$ When $j\geq 1$ the curves belong to the subcase
(ii) considered below. \medskip

\textit{Subcase} (ii$):$ $n=p_{1}d^{\prime }$. Here $\mathcal{E}%
=p_{1}d^{\prime }(q_{1}d^{\prime }-1-k)$ and 
\begin{equation*}
\Delta =(kp_{1}-q_{1}-1)d^{\prime }+2,
\end{equation*}%
which can be non-positive. Of curse, using the changes $\psi \rightarrow
\psi /\varphi ^{j}$ we reduce the problem to the case $m=\mathrm{ord}\psi
_{t=1}<n.$ With $m<n$ we have 
\begin{eqnarray*}
\mathcal{E} &=&(m-1)(n-1)+(n^{\prime }-1)+n^{\prime }(q_{1}d^{\prime }-m) \\
&\leq &(p_{1}d^{\prime }-2)(p_{1}d^{\prime }-1)+(q_{1}d^{\prime
}-p_{1}d^{\prime }+1)n^{\prime }-1
\end{eqnarray*}%
and%
\begin{equation*}
\Delta =(q_{1}d^{\prime }-p_{1}d^{\prime }+1)(p_{1}d^{\prime }-n^{\prime
}-1)-d^{\prime }+2>0
\end{equation*}%
(as $n^{\prime }\leq \frac{1}{2}p_{1}d^{\prime }).$

\bigskip

With this we have completed the proof of Main Theorem for curves of Type $%
\binom{+}{+}.$ $\qed$

\bigskip

\section{Curves of Type $\binom{-+}{+-}$}

Recall that $0<q<p,$ $0<r<s,$ $p+r\leq q+s$ and 
\begin{equation*}
\sigma =\dim Curv=p+r+q+s-3.
\end{equation*}%
Curves of this type are relatively simple due to the fact that $ps-rq$ is
large, it is $\geq 3p^{\prime }r^{\prime }.$

\bigskip

\textbf{6.I. Double points hidden at infinity.} We begin with the situation%
\begin{equation*}
\max n_{i}\leq \max (p^{\prime },r^{\prime }).
\end{equation*}%
Suppose that $p^{\prime }$ is dominating; the case $p^{\prime }<r^{\prime }$
is treated analogously. Then $\mathcal{E}=p^{\prime }(p+r+q+s-3)$ and%
\begin{eqnarray*}
\Delta &=&(pq+2ps+rs-p-r-q-s-p^{\prime }-r^{\prime }+2)-\mathcal{E} \\
&=&(p^{\prime })^{2}(p_{1}q_{1}-p_{1}-q_{1})+p^{\prime }r^{\prime
}(2p_{1}s_{1}-r_{1}-s_{1})+(r^{\prime })^{2}r_{1}s_{1} \\
&&-p^{\prime }(p_{1}+q_{1}-2)-r^{\prime }(r_{1}+s_{1}+1)+2.
\end{eqnarray*}%
If $q_{1}\geq 2,$ than $p_{1}\geq 3$ and $p_{1}q_{1}-p_{1}-q_{1}>0;$ here it
is easy to see that $\Delta >0.$

If $q_{1}=1,$ then we must use the handsomeness property (Definition 2.43
and Proposition 2.44), which ensures that $p^{\prime }<s=s_{1}r^{\prime }$.
Therefore, with $r^{\prime }(r_{1}+1)\leq s,$ we get 
\begin{eqnarray*}
\Delta &\geq &-(p^{\prime })^{2}+p^{\prime }s(2p_{1}-1)-p^{\prime }r^{\prime
}r_{1}+r^{\prime }r_{1}s-p^{\prime }(p_{1}-1)-2s+2 \\
&\geq &(p^{\prime })^{2}(2p_{1}-2)-p^{\prime }(p_{1}+1)+2 \\
&\geq &p^{\prime }(3p_{1}-5)+2>0
\end{eqnarray*}%
(where we used $p^{\prime }\geq 2).$

\bigskip

\textbf{6.II. Double points hidden at the singularities.} We know very well
that $\mathcal{E}$ is maximal when there is only one singular point and that 
\begin{equation*}
\mathcal{E}=n(p+r+q+s-1-n).
\end{equation*}%
(Looking at the below estimates one can see that admitting more
singularities leads to $\Delta >0.)$

We consider two cases:%
\begin{equation*}
\text{A. }n\leq p+r-1,\;\;\;\text{B. }n=p+r.
\end{equation*}
\medskip

\textbf{Case A.} We put $n=p+r-1;$ we shall see that the case is strict and
hence this assumption is correct.

We have $\mathcal{E}=(p+r-1)(q+s)$ and 
\begin{eqnarray*}
\Delta &=&ps-rq-(p+r)-(p^{\prime }+r^{\prime })+2 \\
&=&(p-q-p^{\prime })(s-1)+(s-r-r^{\prime })(q+1)+(p^{\prime
}-1)(s-2)+(r^{\prime }-1)q.
\end{eqnarray*}%
Therefore $\Delta \leq 0$ iff $p=q+p^{\prime },$ $s=r+r^{\prime },$ $%
r^{\prime }=1$ and one of the two:%
\begin{equation*}
\text{A.1. }s=2,\;\;\;\;\;\;\text{A.2. }p^{\prime }=1.
\end{equation*}%
Moreover, the case is strict.

\bigskip

\textbf{A.1.} Here $r=1.$ Thus $p+r=q+p^{\prime }+1$ and $q+s=q+2.$ Since $%
p+r\leq q+s,$ we get $p^{\prime }=1,$ i.e. the case A.2.

\bigskip

\textbf{A.2.} We have $p+r=q+s=n+1$ and $2\delta _{\max }=n^{2}+n$. Since
the case is strict, we should have $\mu =n\nu =\mu _{\min }+n^{\prime }\nu
^{\prime }.$ By Lemma 2.27 either $m=ln+n^{\prime },$ $\nu ^{\prime }=0$ or $%
m=ln.$

If $m<n,$ then $m=n^{\prime }\leq \frac{1}{2}n$ and $\mu =\mu _{\min
}=n(m-1)<2\delta _{\max }.$

Let $m=n,$ i.e. we have the 2-parameter family%
\begin{equation*}
\varphi =(t-1)^{n}(t+\alpha )t^{-r},\;\;\;\psi =(t-1)^{n}(t+\beta )t^{-r-1}.
\end{equation*}%
If $\nu ^{\prime }=0$, i.e. $C_{1}^{(1)}\neq 0,$ then $\mu =n(n-1)<2\delta
_{\max }.$

If $\nu ^{\prime }\geq 1,$ i.e. $C_{1}^{(1)}=0,$ then we make the change $%
\psi \rightarrow \psi -\frac{1+\beta }{1+\alpha }\varphi =-\frac{1+\beta }{%
1+\alpha }(t-1)^{n+2}t^{-r-1}$ and we see that $C_{2}^{(1)}=\frac{1+\beta }{%
1+\alpha }\neq 0.$ If $n$ is odd, then $\mu =n^{2}-1<2\delta _{\max }.$ If $%
n^{\prime }=\gcd (n,n+2)=2,$ then either $C_{3}^{(1)}\neq 0$ (and $\mu
=(n+1)(n-2)+(n+2)(2-1)=n^{2}<2\delta _{\max })$ or $C_{3}^{(1)}=0\neq
C_{5}^{(1)}.$

In the latter case $\mu =n^{2}+2,$ what equals $2\delta _{\max }$ only for $%
n=2.$ Since $n=p+r-1$ and $p\geq 2,$ we get $p=s=2$ and $r=q=1,$ i.e. 
\begin{equation*}
\varphi =(t-1)^{2}(t+\alpha )t^{-1},\;\;\;\psi =(t-1)^{2}(t+\beta )t^{-2}.
\end{equation*}
\medskip

6.1. \textbf{Lemma.} \emph{The above curve is an embedded annulus only for }$%
\alpha =2,$\emph{\ }$\beta =1$\emph{\ (item (w) of Main Theorem).} \medskip

\textit{Proof}. Since $n=2$, there can be other singular points with $%
n_{j}=2.$ There can be either three $\mathbf{A}_{2}$ singularities, or an $%
\mathbf{A}_{4}$ singularity and an $\mathbf{A}_{2}$ singularity, or an $%
\mathbf{A}_{6}$ singularity. With $\alpha ^{\prime }=1/(1+\alpha ),$ $\beta
^{\prime }=1/(1+\beta )$ and $\tau =t-1$ we have 
\begin{eqnarray*}
\alpha ^{\prime }\varphi &=&\tau ^{2}\frac{1+\alpha ^{\prime }\tau }{1+\tau }%
=\tau ^{2}\left[ 1+(\alpha ^{\prime }-1)\tau +\ldots \right] , \\
\beta ^{\prime }\psi &=&\tau ^{2}\frac{1+\beta ^{\prime }\tau }{(1+\tau )^{2}%
}=\tau ^{2}\left[ 1+(\beta ^{\prime }-2)\tau +\ldots \right] .
\end{eqnarray*}%
Here $\frac{d}{d\tau }(\alpha ^{\prime }\varphi )=\tau \left[ 2\alpha
^{\prime }\tau ^{2}+(3\alpha ^{\prime }+1)\tau +2\right] /(1+\tau )^{2}$ and 
$\frac{d}{d\tau }(\beta ^{\prime }\psi )=\tau \lbrack \beta ^{\prime }\tau
^{2}+3\beta ^{\prime }\tau +2]/(1+\tau )^{3}.$

In the case $3\mathbf{A}_{2}$ it should be $2\alpha ^{\prime }\tau
^{2}+(3\alpha ^{\prime }+1)\tau +2\equiv \beta ^{\prime }\tau ^{2}+3\beta
^{\prime }\tau +2.$ This gives $\alpha ^{\prime }=\frac{1}{3}$ and $\beta
^{\prime }=\frac{2}{3},$ i.e. $\alpha =2$ and $\beta =\frac{1}{2};$ it is
the Rudolph's curve [Ru].

Let $\chi :=\alpha ^{\prime }\varphi -\beta ^{\prime }\psi =\tau ^{2}\left[
(1+\alpha ^{\prime }\tau )(1+\tau )-\left( 1+\beta ^{\prime }\tau \right) %
\right] (1+\tau )^{-2}.$ If $C_{1}^{(1)}=0,$ then $\beta ^{\prime }=\alpha
^{\prime }+1$ and $\chi =\alpha ^{\prime }\tau ^{4}(1+\tau )^{-2}=\alpha
^{\prime }\tau ^{4}(1-2\tau +\ldots ).$ Moreover, $\frac{d\chi }{d\tau }%
=2\alpha ^{\prime }\tau ^{3}(2+\tau )(1+\tau )^{-3}.$ If there is another
singular point $t_{2}=1+\tau _{2},$ then $\tau _{2}=-2$ and the condition $%
\frac{d\varphi }{d\tau }(\tau _{2})=0$ gives $\alpha ^{\prime }=0$ (bad
solution). So the case $\mathbf{A}_{2}+\mathbf{A}_{4}$ is not realized.

The possibility $C_{3}^{(1)}=0$ means $2(\alpha ^{\prime }-1)=-2.$ Again we
get the bad solution $\alpha ^{\prime }=0.$ The case $\mathbf{A}_{6}$ is
also not realized. $\qed$

\bigskip

\textbf{Case B:} $n=p+r.$ We have 
\begin{eqnarray*}
\Delta &=&ps-rq-(q+s)-(p^{\prime }+r^{\prime }) \\
&=&(p-q-p^{\prime })(r+1)+(s-r-r^{\prime })(p-1)+(p^{\prime }-1)r+(r^{\prime
}-1)(p-2).
\end{eqnarray*}%
So, the case is strict, $p^{\prime }=1,$ $p=q+1,$ $s=r+r^{\prime }$ and one
of the two:%
\begin{equation*}
\text{B.1. }r^{\prime }=1,\;\;\;\;\;\;\text{B.2. }p=2,\text{ }r^{\prime
}\geq 2.
\end{equation*}%
We also assume that $m\leq q+s-1;$ otherwise, we use Lemma 5.15.

\bigskip

\textbf{B.1.} Here $p+r=q+s=n$ and $2\delta _{\max }=(n-1)^{2}+n-1=n^{2}+n$.
Since $m<n$ and the case is strict, it should be $m=n^{\prime }$ and $\mu
=\mu _{\min }=n(m-1)\leq (p+r)(q+s-2)=n^{2}-2n<2\delta _{\max }.$

\bigskip

\textbf{B.2.} Here $p=2,$ $q=1$ and $r^{\prime }\geq 2.$ If $m\geq n,$ then
we apply the change $\psi \rightarrow \psi /\varphi .$ This gives a curve of
Type $\binom{-}{+}$ (because $\psi /\varphi (\infty )=0).$ But such curves
were already classified.

If $m<n,$ then by the strictness we have $m=n^{\prime }$ and $n\nu =\mu
_{\min }=n(m-1).$ Since $\nu =q+s-1$ and $m-1\leq q+s-2$, we get a
contradiction.

\bigskip

This completes the type $\binom{-+}{+-}$ and finishes the proof of Main
Theorem. $\qed$

\bigskip

\textbf{Acknowledgements}. We would like to thank M. Koras and Z. Jelonek
for their comments about our approach to the affine curves.

\end{document}